\documentclass[12pt]{article}
\usepackage{amssymb, amsmath, amsfonts}
\usepackage{mathrsfs}
\usepackage{dsfont}
\newcommand{\p}{{\sf P}}

\newcommand{\e}{{\sf E}}
\newcommand{\ind}{{\mathbb{I}}}

\newcommand{\var}{{\sf var}}
\newcommand{\cov}{{\sf cov}}

\newcommand\norm[1]{\left\lVert#1\right\rVert}

\newtheorem{thm}{Theorem}

\newtheorem{rem}{Remark}

\newtheorem{cor}{Corollary}

\newtheorem{lem}{Lemma}

\numberwithin{equation}{section}  

\begin{document}

\begin{center}
{\LARGE{\bf Statistical estimation of the Kullback-Leibler 

\vspace{0.2cm}
divergence }}
\end{center}

\begin{center}
{\large {\bf Alexander Bulinski\footnote{E-mail: bulinski@mech.math.msu.su}, Denis Dimitrov\footnote{E-mail: den.dimitrov@gmail.com}}}
\end{center}
\begin{center}
{\it Dept. of Mathematics and Mechanics, Lomonosov Moscow State University,\\ Moscow 119234, Russia}
\end{center}
\vskip0.5cm
{\small
{\bf Abstract} Wide conditions are provided to guarantee asymptotic unbiasedness and $L^2$-consistency of the introduced estimates of the Kullback - Leibler divergence for probability measures in $\mathbb{R}^d$ having densities w.r.t. the Lebesgue measure. These estimates are constructed by means of two independent collections of i.i.d. observations and
involve the specified $k$-nearest neighbor statistics. In particular, the established results are valid
for estimates of the Kullback - Leibler divergence between any two Gaussian measures in $\mathbb{R}^d$ with nondegenerate covariance matrices.
As a byproduct we obtain new statements concerning the Kozachenko-Leonenko estimators of the
Shannon differential entropy.}

\vskip0.2cm
\noindent
{\bf AMS (2010) Subject Classification} 60F25, 62G20, 62H12

\section{Introduction}

The Kullback - Leibler divergence
plays important role in various domains such as
statistical inference
(see, e.g., \cite{Moulin}, \cite{Pardo}), machine learning (\cite{Bishop}, \cite{Poczos}), computer vision (\cite{Cui}, \cite{Deledalle}),  network security (\cite{Ma}, \cite{Yu}), feature selection and classification (\cite{Li}, \cite{Peng}, \cite{Vergara}), physics (\cite{Granero}), biology
(\cite{Charzynska}), finance (\cite{Zhou}), among others.
Recall that this divergence measure between probabilities $\mathbb{P}$ and $\mathbb{Q}$ on a space $(S,\mathcal{B})$
is defined by way of
\begin{equation}\label{divergence}
    D(\mathbb{P}||\mathbb{Q}) :=
\int\limits_{\it S} \log\left(\frac{d\mathbb{P}}{d\mathbb{Q}}\right) d\mathbb{P}\;\;\mbox{if}\;\;\mathbb{P}\ll \mathbb{Q},
\end{equation}
where $\frac{d\mathbb{P}}{d\mathbb{Q}}$ stands for the Radon-Nikodym derivative. Otherwise, $D(\mathbb{P}||\mathbb{Q}) := +\infty$.
We employ the base $e$ of logarithms (a constant factor is not essential here).
It is worth to emphasize that mutual information, widely used in many research directions, is a special case of
the Kullback -Leibler divergence for certain measures.
For comparison of various $f$-divergence measures see \cite{Sason}.

If $(S,\mathcal{B}) = (\mathbb{R}^d, \mathcal{B}(\mathbb{R}^d))$ and (absolutely continuous) $\mathbb{P}$ and $\mathbb{Q}$ have densities, $p(x)$ and $q(x)$, $x\in \mathbb{R}^d$, w.r.t. the Lebesgue measure $\mu$, then \eqref{divergence} can be rewritten as
\begin{equation}\label{divergence_densities}
    D({\mathbb{P}}||{\mathbb{Q}}) = \int\limits_{\mathbb{R}^d}p(x)\log\left(\frac{p(x)}{q(x)}\right) dx \;\; \mbox{for}\;\;\mathbb{P} \ll \mathbb{Q},
\end{equation}
otherwise, $D({\mathbb{P}}||{\mathbb{Q}}) = +\infty$.
To simplify notation we write $dx$ instead of $\mu(dx)$. We formally set
$0/0:=0$,
$0\cdot\log 0:=0$.
For a (version of) probability density $f$ denote by
$S(f):=\{x\in \mathbb{R}^d: f(x)>0\}$  its support. Clearly,  the integral in \eqref{divergence_densities}
is taken over $S(p)$.
Observe that when $\mathbb{P} \ll \mu$ and $\mathbb{Q} \ll \mu$ then $\mathbb{P} \ll \mathbb{Q}$ if and only if $\mathbb{P}(S(p) \setminus S(q)) = 0$.
Formula \eqref{divergence_densities} is closely related to
cross-entropy and the Shannon differential entropy.

Usually one  has to reconstruct the measures (describing a stochastic model under consideration) or their characteristics using some collections of observations. In the pioneering paper \cite{Kozachenko} the estimator of the Shannon differential entropy was proposed, based on the nearest neighbor statistics. In a series of papers this estimate was
studied and applied. Moreover, estimators of the R\'enyi entropy, mutual information and the Kullback - Leibler divergence have appeared (see, e.g., \cite{Kraskov}, \cite{Leonenko},  \cite{Wang}).
However, the authors of \cite{Pal} indicated the occurrence of gaps in the known proofs concerning the limit behavior of such statistics. This issue has attracted our attention and motivated our study of the declared asymptotic properties. Thus in a recent work \cite{Bul_Dim} the new functionals were introduced to prove asymptotic unbiasedness and $L^2$-consistency of the Kozachenko - Leonenko estimators of the Shannon differential entropy.
The present paper is aimed at extension of our approach to grasp the Kullback - Leibler divergence estimation. Instead of the nearest neighbor statistics  we employ the $k$-nearest neighbor statistics (on order statistics see, e.g., \cite{Biau}) and also use more general forms of the mentioned functionals.

Let $X$ and $Y$ be random vectors taking values in $\mathbb{R}^d$ and having distributions $\p_X$ and $\p_Y$, respectively (further we consider $\mathbb{P}=\p_X$ and $\mathbb{Q}=\p_Y$). Consider i.i.d. random vectors $X_1,X_2,\ldots,$ and i.i.d. random vectors $Y_1,Y_2,\ldots,$ with $law(X_1) = law(X)$ and $law(Y_1)=law(Y)$. Assume that $\{X_i, Y_i, i\in \mathbb{N}\}$ are independent. We are interested in statistical estimation
of $D(\p_X||\p_Y)$ constructed by means of  observations $\mathbb{X}_n:=\{X_1,\ldots,X_n\}$ and $\mathbb{Y}_m:=\{Y_1,\ldots,Y_m\}$,  $n,m\in \mathbb{N}$.
All random variables under consideration are defined on a complete probability space $(\Omega,\mathcal{F},\p)$.

For a finite set $E=\{z_1,\ldots,z_N\}\subset \mathbb{R}^d$, where $z_i\neq z_j$ $(i\neq j)$,
and a vector $v\in \mathbb{R}^d$, renumerate points of $E$ as $z_{(1)}(v),\ldots,z_{(N)}(v)$ in such a way that $\|v-z_{(1)}\|\leq \ldots \leq \|v-z_{(N)}\|$, here $\|\cdot\|$ is the Euclidean norm in
$\mathbb{R}^d$. If there are points $z_{i_1},\ldots,z_{i_s}$ having the same distance from $v$
then we numerate them according the increasing indexes among $i_1,\ldots,i_s$.
In other words, for $k=1,\ldots,N$,  $z_{(k)}(v)$ is the $k$-NN (Nearest Neighbor) for $v$
in a set $E$. To indicate that $z_{(k)}(v)$ is constructed by means of $E$ we write $z_{(k)}(v,E)$.
Fix $k\in \{1,\ldots,n-1\}$, $l\in \{1,\ldots,m\}$ and (for each $\omega \in \Omega$) put
$$
R_{n,k}(i):= \|X_i-X_{(k)}(X_i,\mathbb{X}_n\setminus \{X_i\})\|,\;\;\;
V_{m,l}(i):= \|X_i-Y_{(l)}(X_i,\mathbb{Y}_m)\|,\;\;i=1,\ldots,n.
$$
We assume that $X$ and $Y$ have  densities $p=\frac{d\p_X}{d\mu}$ and $q=\frac{d\p_Y}{d\mu}$. Then with probability one all points in $\mathbb{X}_n$ are distinct as well as points of $\mathbb{Y}_m$.

Introduce an estimate of $D(\p_X||\p_Y)$, for $n\geq k+1$ and $m\geq l$, letting
\begin{equation}\label{kl_div_estimation}
\widehat{D}_{n,m}(k,l) := \psi{(k)}-\psi{(l)} + \frac{1}{n} \sum_{i=1}^n \log\left(\frac{m V^d_{m,l}(i)}{(n-1) R^d_{n,k}(i)}\right).
\end{equation}
Here $\psi(t) = \frac{d}{dt} \log{\Gamma(t)} = \frac{\Gamma'(t)}{\Gamma(t)}$ is the digamma function, $t>0$.

\begin{rem}\label{r1}
\normalfont{ If $k=l$ then
$$
\widehat{D}_{n,m}(k) = \frac{d}{n}\sum_{i=1}^n \log\left(\frac{V_{m,l}(i)}{R_{n,k}(i)}\right)
+ \log\left(\frac{m}{n-1}\right),
$$
and we come to formula (5) in \cite{Wang}}.
\end{rem}
\vspace{-0.5cm}
\begin{rem}\label{remark_generalize}
\normalfont{ All our results will be valid for the following generalization of statistics $\widehat{D}_{n,m}(k,l)$:
\begin{equation}\label{kl_div_estimation_most_general_case_with_ki}
    \widetilde{D}_{n,m}(\mathcal{K}_n,\mathcal{L}_n) := \frac{1}{n} \sum_{i=1}^n \left(\psi{(k_i)}-\psi{(l_i)}\right) + \log\left(\frac{m}{n-1}\right) \\ + \frac{d}{n} \sum_{i=1}^n \log\left(\frac{V_{m,l_i}(i)}{R_{n,k_i}(i)}\right),
\end{equation}
where $\mathcal{K}_n := \{k_i\}_{i=1}^n$, $\mathcal{L}_n := \{l_i\}_{i=1}^n$ and, for some $r \in \mathbb{N}$ and all $i \in \mathbb{N}$, $k_i \leq r$, $l_i \leq r$. Note that \eqref{kl_div_estimation_most_general_case_with_ki} is well-defined for $n \geq \max_{i=1,\ldots,n} k_i + 1$, $m \geq \max_{i=1,\ldots,n}l_i$. We will only  consider  the estimates \eqref{kl_div_estimation} since the study of $\widetilde{D}_{n,m}(\mathcal{K}_n,\mathcal{L}_n)$ follows the same lines.}
\end{rem}

\vspace{-0.3cm}
Developing the approach of \cite{Bul_Dim} to analysis of asymptotic behavior of the Kozachenko-Leonenko estimates of the Shannon differential entropy (introduced in \cite{Shannon}, Part III, Section 20) we encounter new complications due to dealing with $k$-nearest neighbor statistics for $k\in \mathbb{N}$ (not only for
$k=1$). Accordingly, in the framework of the Kullback-Leibler divergence estimation, we propose a new way to bound the function $1- F_{m,l,x}(u)$ playing the key role in the proofs (see formula \eqref{eq1a}). Also instead of the function $G(t)=t\log t$ (for $t>1$), used in \cite{Bul_Dim} for study of the Shannon entropy estimates, we employ a regularly varying function
$G_N(t)= t \log_{[N]}(t)$ where (for $t$ large enough) $\log_{[N]}(t)$ is the $N$-fold iteration of the logarithmic function and $N\in \mathbb{N}$ is chosen arbitrarily. Whence in
the definition of integral functional $K_{p,q}(\nu,N,t)$ by formula \eqref{p1} below one can take a function $G_N(z)$ having, for $z>0$, the growth rate close to that of function $z$.  Moreover, this permits a generali\-zation
of \cite{Bul_Dim} results. Here we invoke convexity of $G_N$
(see Lemma~\ref{l6}) to provide more simple conditions for asymptotic unbiasedness and
$L^2$-consistency of the Shannon differential entropy
 than those employed in \cite{Bul_Dim}.

Mention in passing that there exist  investigations treating other important aspects of the
mutual information and entropy estimation.
In \cite{Alonso-Ruiz} entropy estimators are applied to detection of the fiber materials inhomogeneities. The mixed models and conditional entropy estimation are studied, e.g., in
\cite{Bul_Koz}, \cite{Coelho}. The central limit theorem for the Kozachenko-Leonenko estimates is established in \cite{Delattre}. The limit theorems for point processes on manifolds are employed
in \cite{Penrose} to analyze behavior of the Shannon and the R\'enyi entropy estimates.
The convergence rates for the Shannon entropy (truncated) estimates are obtained in
\cite{Tsybakov} for one-dimensional case, see also \cite{Singh_1} for multidimensional case.
Ensemble estimation of density functional is considered in \cite{Sricharan}.
A recursive rectilinear partitioning for the differential entropy is considered in
\cite{Stowell}. The mutual information estimation by the local Gaussian approximation is developed in \cite{Gao}.
Note that various deep results (including the central limit theorem) were obtained for the Kullback - Leibler estimates under certain conditions imposed on derivatives of
unknown densities (see, e.g., the  recent papers \cite{Berrett}, \cite{Moon}, \cite{Sasaki}). Our goal is to provide wide conditions for the asymptotic unbiasedness and $L^2$-consistency of the Kullback - Leibler divergence estimates \eqref{kl_div_estimation}, as $n,m\to \infty$, without such smoothness hypothesis. Also we do not assume that densities have bounded supports.

The paper is organized as follows. In Section 2 we formulate main results, Theorems \ref{th1} and \ref{th_main2}. Their proofs are presented in
Sections 3 and 4, respectively. Proofs of several lemmas are given in Appendix (Section 5).

\section{Main results}

Some notation is necessary.
For a probability density $f$ in $\mathbb{R}^d$,
$x\in \mathbb{R}^d$, $r>0$ and $R>0$, as in \cite{Bul_Dim}, introduce the functions (or functionals depending on parameters)
\begin{equation}\label{I}
I_f(x,r):= \frac{\int_{B(x,r)} f(y) \, dy}{r^d V_d},
\end{equation}
\begin{equation}\label{mM}
M_f(x, R) := \sup_{r \in (0,R]} I_f(x,r),\;\;m_f(x, R) := \inf_{r \in (0,R]} I_f(x,r),
\end{equation}
where $B(x,r):=\{y\in \mathbb{R}^d: \|x-y\|\leq r\}$.
Observe that changing $\sup_{r \in (0,R]}$ by $\sup_{r \in (0,\infty)}$
in the definition of $M_f(x,R)$ leads to the celebrated Hardy - Littlewood maximal
function $M_f(x)$ widely used in harmonic analysis.
Some properties
of the function $\int_{B(x,r)} f(y) \, dy$ are considered, e.g., in \cite{Evans}.
According to Lemma 2.1 \cite{Bul_Dim}, for a probability density $f$ in $\mathbb{R}^d$, the function $I_f(x,r)$ defined in \eqref{I} is continuous in $(x,r)\in \mathbb{R}^d\times (0,\infty)$.

Set $e_{[0]}:=1$ and  $e_{[N]}:=\exp\{e_{[N-1]}\}$, $N\in \mathbb{N}$.
Introduce a function $\log_{[1]}(t):=\log t$, $t>0$.
For $N\in \mathbb{N}$, $N> 1$, set
$\log_{[N]}(t):= \log(\log_{[N-1]}(t)).
$
Evidently, this function (for $N>1$) is defined if  $t > e_{[N-2]}$.
For $N\in \mathbb{N}$, consider the continuous nondecreasing function $G_N: \mathbb{R}_+\to \mathbb{R}_+$, given by formula
\begin{equation}\label{GN}
G_N(t):=\begin{cases}
0, &t\in [0,e_{[N-1]}],\\
t \log_{[N]}(t), &t \in (e_{[N-1]}, \infty).
\end{cases}
\end{equation}

For probability densities $p, q$ in $\mathbb{R}^d$, some $N \in \mathbb{N}$ and positive constants $\nu, t, \varepsilon, R$,  we define the following functionals
with values in $[0,\infty]$
\begin{equation}\label{p1}
K_{p,q}(\nu, N, t):=\;\;\;\;\;\;\;\int\!\!\!\!\!\!\!\!\!\!\!\!\!\!\!\!\!\int\limits_{x,y\in \mathbb{R}^d\!,\,\|x-y\|> t} {G_N\big(|\log \norm{x-y}|^{\nu}\big)} p(x) q(y) \, dx \, dy,
\end{equation}
\begin{equation}\label{p2}
Q_{p,q}(\varepsilon, R):=\int_{\mathbb{R}^d} M_q^{\varepsilon}(x, R) p(x)\,dx,
\end{equation}
\begin{equation}\label{p3}
T_{p,q}(\varepsilon,R):=\int_{\mathbb{R}^d} m_q^{-\varepsilon}(x, R) p(x) \, dx.
\end{equation}
Set $K_{p,q}(\nu, N)  := K_{p,q}(\nu, N, e_{[N]})$. Clearly, for any $N\in\mathbb{N}$, $\nu,t,u>0$
such that $t<u$, one has
\begin{equation}\label{K}
K_{p,q}(\nu, N, u)\leq K_{p,q}(\nu, N, t)\leq K_{p,q}(\nu, N, u) + \max\{G_N(|\log t|^{\nu}), G_N(|\log u|^{\nu})\}.
\end{equation}

\begin{rem}\label{rem1}
\normalfont{
We stipulate that
$1/0:=\infty$ (consequently $m_q^{-\varepsilon_2}(x, R):=\infty$ when $m_q(x,R)=0$).
For arbitrary versions of  $p$ and $q$, we can write in
\eqref{p2}, \eqref{p3}  the integrals over the support $S(p)$ instead of integrating over $\mathbb{R}^d$ (obviously, the results do not depend on the choice of  versions).}
\end{rem}

\begin{thm}\label{th1}
Let $\p_X$ and $\p_Y$ have densities $p$ and $q$, respectively.
Suppose that $p$ and $q$ are such that, for some $\varepsilon_i>0, R_i>0$ and  $N_j \in \mathbb{N}$, where $i = 1,2,3,4$ and $j=1,2$, the functionals $K_{p,q}(1, N_1)$, $Q_{p,q}(\varepsilon_1, R_1)$, $T_{p,q}(\varepsilon_2, R_2)$, $K_{p,p}(1, N_2)$, $Q_{p,p}(\varepsilon_3, R_3)$, $T_{p,p}(\varepsilon_4, R_4)$ are finite.
Then, for any fixed  $k,l\in \mathbb{N}$, the estimates $\widehat{D}_{n,m}(k,l)$, introduced in \eqref{kl_div_estimation},  are asymptotically unbiased, i.e.
\begin{equation}\label{main1}
\lim\limits_{n,m\to \infty}\e \widehat{D}_{n,m}(k,l) = D(\p_X || \p_Y).
\end{equation}
\end{thm}

\begin{rem}\label{rem1add}
It is useful to note that if $Q_{p,q}(\varepsilon_1, R_1) < \infty$ and $T_{p,q}(\varepsilon_2, R_2) < \infty$ for some positive $\varepsilon_1, \, \varepsilon_2, R_1, R_2$ then $\int_{\mathbb{R}^d} p(x)  |\log{q(x)}| \, dx < \infty$.
Indeed, definition \eqref{mM}
and the Lebesgue differentiation theorem (see, e.g., Theorem 25.17 \cite{Yeh})
yield that $m_q(x, R_2) \leq q(x) \leq M_q(x, R_1)$ for $\mu$-almost all $x \in \mathbb{R}^d$.
Evidently, $\log z \leq \frac{1}{\varepsilon}z^{\varepsilon}$ for any $z\geq 1$ and each $\varepsilon >0$. Consequently,
\begin{gather*}\int_{\mathbb{R}^d} p(x)  |\log{q(x)}| \, dx = \int_{q(x) \geq 1} p(x)  \log{q(x)} \, dx + \int_{q(x) < 1} p(x)  \log{\frac{1}{q(x)}} \, dx \\ \leq  \frac{1}{\varepsilon_1} Q_{p,q}(\varepsilon_1, R_1) + \frac{1}{\varepsilon_2} T_{p,q}(\varepsilon_2, R_2) < \infty.
\end{gather*}
So, the integrals $Q_{p,q}(\varepsilon_1, R_1)$, $T_{p,q}(\varepsilon_2, R_2)$, $Q_{p,p}(\varepsilon_3, R_3)$, $T_{p,p}(\varepsilon_4, R_4)$ finiteness implies
the finiteness of integral in \eqref{divergence_densities}
(and also guarantees that $\p_X\ll \p_Y$).
\end{rem}

\begin{lem}\label{lemma1} Let $p$ and $q$ be any probability densities in $\mathbb{R}^d$.
Then the following statements are valid.

\noindent
$1)$ If $K_{p,q}(\nu_0,N_0)<\infty$ for some $\nu_0 > 0$ and $N_0 \in \mathbb{N}$ then $K_{p,q}(\nu,N)<\infty$ for any $\nu \in (0, \nu_0]$ and each $N \geq N_0$.

\noindent
$2)$ If $Q_{p,q}(\varepsilon_1, R_1)<\infty$ for some $\varepsilon_1 > 0$ and $R_1 > 0$ then $Q_{p,q}(\varepsilon, R)<\infty$ for any $\varepsilon \in (0,\varepsilon_1]$ and each $R>0$.

\noindent
$3)$ If $T_{p,q}(\varepsilon_2,R_2)<\infty$ for some $\varepsilon_2>0$ and $R_2>0$ then $T_{p,q}(\varepsilon,R)<\infty$ for any $\varepsilon \in (0,\varepsilon_2]$ and each $R>0$.
\end{lem}

The proof is given in Appendix. In view of Lemma \ref{lemma1}, one can recast Theorem \ref{th1} as follows.

\begin{cor}\label{cor1} Let, for some positive $\varepsilon, R$ and $N \in \mathbb{N}$, the functionals $K_{p,q}(1, N)$, $Q_{p,q}(\varepsilon, R)$, $T_{p,q}(\varepsilon, R)$, $K_{p,p}(1, N)$, $Q_{p,p}(\varepsilon, R)$, $T_{p,p}(\varepsilon, R)$ be finite.
Then \eqref{main1} holds. Moreover, we obtain the equivalent conditions assuming
that these functionals are finite for some $\varepsilon >0$ and $R=\varepsilon$.
\end{cor}

Let us also consider the following simple conditions.
\vspace{0.3cm}

\noindent
$(A;p,q,\nu)$ For  probability densities $p, q$ in $\mathbb{R}^d$ and some positive $\nu$
\begin{equation}\label{pq1}
L_{p,q}(\nu) := \int_{\mathbb{R}^d}\int_{\mathbb{R}^d}|\log \norm{x-y}|^{\nu} p(x)q(y)\,dxdy <\infty.
\end{equation}
We formally set $\log 0:=-\infty$ and, as usual, 
$\int_A g(z) Q(dz)=0$ whenever $g(z)= \infty$ (or $-\infty$) for $z\in A$ and $Q(A)=0$,
where $Q$ is a $\sigma$-finite measure on $(\mathbb{R}^d,\mathcal{B}(\mathbb{R}^d))$.

\vskip0.2cm
\noindent
$(B_1; f)$ There exists a version of density $f$ such that, for some $M(f)\in (0,\infty)$,
\begin{equation*}
f(x)\leq M(f),\;\;x\in \mathbb{R}^d.
\end{equation*}

\noindent
$(C_1; f$) There exists a version of density $f$ such that, for some $m(f)\in (0,\infty)$,
\begin{equation*}
f(x)\geq m(f),\;\;x\in S(f).
\end{equation*}

\begin{cor}\label{cor2}
Let conditions $(A;p,q,\nu)$ and $(A;p,p,\nu)$ be satisfied with some $\nu>1$.
Then \eqref{main1} is true, provided that $(B_1;f)$ and $(C_1;f)$ are valid for $f=p$ and $f=q$.
Moreover, if the latter assumption concerning $(B_1;f)$ and $(C_1;f)$ holds then \eqref{main1} is true whenever $p$ and $q$ have  bounded supports.
\end{cor}

Next we formulate conditions to guarantee $L^2$-consistency of estimates \eqref{kl_div_estimation}.
\begin{thm}\label{th_main2}
Let the requirements $K_{p,q}(1,N_1)<\infty$ and
$K_{p,p}(1,N_2)<\infty$ in conditions of Theorem \ref{th1} be replaced by $K_{p,q}(2,N_1)<\infty$ and
$K_{p,p}(2,N_2)<\infty$.
Then, for any fixed $k,l \in \mathbb{N}$, the estimates $\widehat{D}_{n,m}(k,l)$ are $L^2$-consistent, i.e.
\begin{equation}\label{main2}
\lim\limits_{n,m\to \infty}
\e \left(\widehat{D}_{n,m}(k,l) - D(\p_X || \p_Y)\right)^2 = 0.
\end{equation}
\end{thm}

Due to Lemma \ref{lemma1} one can recast Theorem \ref{th_main2} as follows.

\begin{cor}\label{cor1aa}
Let, for some positive $\varepsilon, R$ and  $N \in \mathbb{N}$, the functionals $K_{p,q}(2, N)$, $Q_{p,q}(\varepsilon, R)$, $T_{p,q}(\varepsilon, R)$, $K_{p,p}(2, N)$, $Q_{p,p}(\varepsilon, R)$, $T_{p,p}(\varepsilon, R)$ be finite.
Then \eqref{main2} holds. Moreover, we obtain the equivalent conditions assuming
that these functionals are finite for some $\varepsilon >0$ and $R=\varepsilon$.
\end{cor}

\begin{cor}\label{cor3}
Let conditions $(A;p,q,\nu)$ and $(A;p,p,\nu)$ be satisfied with some $\nu>2$. Assume that
$(B_1;f)$ and $(C_1;f)$ are valid for $f=p$ and $f=q$.
Then \eqref{main2} is true.
Moreover, if the latter assumption concerning $(B_1;f)$ and $(C_1;f)$ holds then \eqref{main2} is true whenever $p$ and $q$ have  bounded supports.
\end{cor}

\vspace{0.1cm}

Note that D.Evans  considered the ``positive density condition''
in Definition 2.1 of \cite{Evans}
meaning that there exist constants $\beta >1$ and $\delta >0$ such that $\frac{r^d}{\beta}\leq \int_{B(x,r)}q(y)dy \leq \beta r^d $ for all $0\leq r\leq \delta$ and $x\in \mathbb{R}^d$.
Consequently $m_q(x,\delta)\geq \frac{1}{\beta V_d} := m > 0$, $x\in \mathbb{R}^d$. Then $T_{p,q}(\varepsilon, \delta) \leq m^{-\varepsilon} \int_{\mathbb{R}^d} p(x) \, dx = m^{-\varepsilon} < \infty$ for all $\varepsilon > 0$. Analogously, $M_q(x,\delta) \leq \frac{\beta}{V_d} := M$, $M > 0$, $x \in \mathbb{R}^d$, and $Q_{p,q}(\varepsilon, \delta) \leq M^{\varepsilon} \int_{\mathbb{R}^d} p(x) \, dx = M^{\varepsilon} < \infty$ for all $\varepsilon > 0$.
It was proved in \cite{Evans_1} that if $f$ is smooth and its support is a compact convex body in $\mathbb{R}^d$ then the mentioned inequalities from Definition 2.1 of \cite{Evans} hold. Therefore, if $p$ and $q$ are smooth and their supports are  compact convex bodies in $\mathbb{R}^d$ then 
one can simplify conditions of Corollaries \ref{cor1} and \ref{cor1aa}.

\vspace{0.3cm}
Now instead of (C1; $f$) we consider the following condition introduced in \cite{Bul_Dim} that allows us to work with densities, whose supports need not be bounded.

\vspace{0.1cm}
\noindent
$(C_2; f)$ For a fixed $R>0$, there exists a constant $c>0$ and a version of a density $f$ such that
\begin{equation}\label{reg}
m_f(x,R)\geq c f(x),\;\;x\in \mathbb{R}^d.
\end{equation}

\begin{rem}\label{rem2a}
\normalfont{
If, for some positive $\varepsilon$, $R$ and $c$, condition $(C_2; q$) is true and
\begin{equation}\label{densi}
\int_{\mathbb{R}^d} q(x)^{-\varepsilon} p(x) dx <\infty,
\end{equation}
then obviously $T_{p,q}(\varepsilon,R)<\infty$.
Thus in Theorems \ref{th1} and \ref{th_main2} one can employ, for $f=p$ and $f=q$, condition $(C_2;f)$ and suppose, for some $\varepsilon >0$, finiteness of $\int_{\mathbb{R}^d} q(x)^{-\varepsilon} p(x) dx$ and $\int_{\mathbb{R}^d} p^{1-\varepsilon}(x) dx$ instead of the corresponding assumptions
 $T_{p,q}(\varepsilon,R)<\infty$ and $T_{p,p}(\varepsilon, R)<\infty$. To illustrate this observation
we provide a result for a density with unbounded support.}
\end{rem}

\begin{cor}\label{cor4}
Let $X$, $Y$ be Gaussian random vectors in $\mathbb{R}^d$ with $\e X = \mu_X$, $\e Y = \mu_Y$
and  nondegenerate covariance matrices $\Sigma_X$ and $\Sigma_Y$, respectively.
Then relations \eqref{main1} and \eqref{main2} hold where
$$
D(\p_X || \p_Y) = \frac{1}{2} \left( {\sf tr}\left( \Sigma_Y^{-1} \Sigma_X \right) + \left( \mu_Y - \mu_X \right)^T \Sigma_Y^{-1} \left( \mu_Y - \mu_X \right) - d + \log{\left( \frac{\det\Sigma_Y}{\det\Sigma_X} \right)} \right).
$$
\end{cor}
The latter formula can be found, e.g., in  \cite{Moulin}, p. 147. The proof of Corollary \ref{cor4} is discussed in Appendix.

Similarly to condition $(C_2;f)$ let us consider the following one.

\vspace{0.3cm}
\noindent
$(B_2; f)$ For a fixed $R>0$, there exists a constant $C>0$ and a version of a density $f$ such that
\begin{equation}\label{reg_M}
M_f(x,R)\leq C f(x),\;\;x\in S(f).
\end{equation}

\begin{rem}\label{rem2aaa}
\normalfont{
If, for some positive $\varepsilon$, $R$ and $c$, condition $(B_2;q)$ is true and
\begin{equation}\label{densi2aaa}
\int_{\mathbb{R}^d} q(x)^{\varepsilon} p(x) dx <\infty
\end{equation}
then obviously $Q_{p,q}(\varepsilon,R)<\infty$.
Thus in Theorems \ref{th1} and \ref{th_main2} one can employ, for $f=p$ and $f=q$, condition $(B_2;f)$ and suppose that $\int_{\mathbb{R}^d} q(x)^{\varepsilon} p(x) dx$ and $\int_{\mathbb{R}^d} p^{1+\varepsilon}(x) dx$ are finite (for some $\varepsilon >0$) instead of the assumptions
$Q_{p,q}(\varepsilon,R)<\infty$ and $Q_{p,p}(\varepsilon, R)<\infty$.}
\end{rem}

For a fixed $k\in \{1,\ldots,n-1\}$, consider the Kozachenko - Leonenko estimate of the Shannon differential entropy $H(X)$ of a vector $X$ with values in $\mathbb{R}^d$ having a density $p$ w.r.t. the Lebesgue measure.
Namely, $H(X):= -\int_{\mathbb{R}^d} (\log p(x)) p(x)\mu(dx)$ and, for i.i.d. observations $X_1,X_2,\ldots$, such that $law (X_1)=law(X)$, set for all $n \geq k+1$,
\vspace{-0.2cm}
\begin{equation}\label{entropy_estimation_general_k}
\widehat{H}_n(k):=\frac{1}{n}\sum_{i=1}^n \log \left(\frac{R_{n,k}^d(i) V_d (n-1)}{e^{\psi(k)}}\right).
\end{equation}
Similar to \eqref{kl_div_estimation_most_general_case_with_ki} one can employ the following generalization of statistics $\widehat{H}_{n}(k)$:
\begin{equation*}
\widetilde{H}_n(\mathcal{K}_n):=-\frac{1}{n} \sum_{i=1}^n \psi{(k_i)} + \log{V_d} + \log{(n-1)} + \frac{d}{n}\sum_{i=1}^n \log R_{n,k_i}(i),
\end{equation*}
where $\mathcal{K}_n := \{k_i\}_{i=1}^n$, and, for some $r \in \mathbb{N}$ and all $i \in \mathbb{N}$, $k_i \leq r$.

\begin{cor}\label{ShEn}
Let  $Q_{p,p}(\varepsilon,R)<\infty$ and $T_{p,p}(\varepsilon,R)<\infty$ for some positive $\varepsilon$ and $R$. Then the following statements hold for any fixed $k \in \mathbb{N}$.

\noindent
1) If, for some $N \in \mathbb{N}$, $K_{p,p}(1, N)<\infty$, then \,
$
\e \widehat{H}_n(k)\to H(X),\;\;n\to \infty.
$

\noindent
2) If, for some $N \in \mathbb{N}$, $K_{p,p}(2, N)<\infty$, then \,
$
\e (\widehat{H}_n(k) - H(X))^2 \to 0,\;\;n\to \infty.
$

In particular, one can employ $L_{p,p}(\nu)$ with $\nu > 1$ instead of $K(1, N)$, and with  $\nu > 2$  instead of $K(2, N)$, where $N \in \mathbb{N}$.
\end{cor}

The proof of the first statement of this corollary is contained in the proof of Theorem \ref{th1}, \textit{Step 5}. In a similar way one can infer the second statement of Corollary \ref{ShEn} by means of the proof of Theorem \ref{th_main2},~\textit{Step~5}.

\section{Proof of Theorem \ref{th1}}

For $n, m \in \mathbb{N}$ such that $n>1$, for  fixed $k\in \mathbb{N}$ and $m \in \mathbb{N}$, where $1 \leq k \leq n-1$, $1 \leq l \leq m$ and
$i = 1, \ldots, n$, set
$\phi_{m,l}(i) = m V^d_{m,l}(i)$, $\zeta_{n,k}(i) = (n-1) R^d_{n,k}(i)$.
Then we can rewrite the estimate $\widehat{D}_{n,m}(k,l)$ as follows
\begin{equation}\label{a1}
\widehat{D}_{n,m}(k,l) = \psi(k) - \psi(l) + \frac{1}{n} \sum_{i=1}^n \big( \log\phi_{m,l}(i) - \log\zeta_{n,k}(i)\big).
\end{equation}

It is sufficient to prove the following two claims.

{\it Statement 1}. For each fixed $l$, all $m$ large enough and any $i \in \mathbb{N}$,
$\e |\log\phi_{m,l}(i)|$ is finite. Moreover,
\begin{equation}\label{mr}
\frac{1}{n} \sum_{i=1}^n \log\phi_{m,l}(i)=\e \log\phi_{m,l}(1) \to \psi(l) - \log{V_d} - \int_{\mathbb{R}^d} p(x) \log{q(x)} \, dx, \;\; m \to \infty.
\end{equation}

{\it Statement 2}. For each fixed $k$, all $n$ large enough and any $i\in \mathbb{N}$,
$\e |\log\zeta_{n,k}(i)|$ is finite. Moreover,
\begin{equation}\label{mra}
\frac{1}{n} \sum_{i=1}^n \log\zeta_{m,l}(i)=\e \log\zeta_{n,k}(1) \to \psi(k) - \log{V_d} - \int_{\mathbb{R}^d} p(x) \log{p(x)} \, dx, \;\; n \to \infty.
\end{equation}

Then in view of \eqref{a1}, \eqref{mr} and \eqref{mra}
\begin{gather*}
\e \widehat{D}_{n,m}(k,l) \to - \int_{\mathbb{R}^d} p(x) \log{q(x)} \, dx +  \int_{\mathbb{R}^d} p(x) \log{p(x)} \, dx
= D(\p_X || \p_Y), \;\; n,m \to \infty.
\end{gather*}

We are going to discuss in detail only the proof of \textit{Statement 1}, since  \textit{Statement 2} is established in a similar way.
It was explained in \cite{Bul_Dim} that if $V$ is a nonegative random variable (hence $\e V \leq \infty$) and $X$ is an arbitrary random vector with values in $\mathbb{R}^d$ then
\begin{equation}\label{ce}
\e V=\int_{\mathbb{R}^d} \e (V|X=x)\p_X(dx).
\end{equation}
Formula \eqref{ce} means that
simultaneously both sides are finite or infinite and coincide.
Let $F(u,\omega)$ be a regular conditional distribution function
of $V$ given $X$ where $u \in [0,\infty)$ and $\omega \in \Omega$.
Let $h$ be a measurable function such that $h:\mathbb{R}\to [0,\infty)$. Then, for $\p_{X}$-almost all $x\in \mathbb{R}^d$, it follows (without assumption $\e h(V)<\infty$) that
\begin{equation}\label{ce1}
\e (h(V)|X=x)= \int_{[0,\infty)}h(u)dF(u,x).
\end{equation}
This means that both sides of \eqref{ce1} are finite or infinite simultaneously and coincide.

By virtue of \eqref{ce} and \eqref{ce1} one can prove that $\e |\log\phi_{m,l}(i)| < \infty$, for all $m$ large enough,
fixed $l$ and for all $i\in \mathbb{N}$,
and \eqref{mr} holds. For this purpose we take $V=\phi_{m,l}(i)$,
$X=X_i$ and $h(u)= |\log u|$, $u>0$ (we use $h(u)=\log^2 u$ in the proof of Theorem \ref{th_main2}).
To reduce the volume of the paper we
only consider below the evaluation of $\e \log\phi_{m,l}(i)$ as all steps of the proof are the same when
treating $\e|\log\phi_{m,l}(i)|$.

We divide the proof of  Statement 1  into four steps. Preliminary \textit{Steps 1-3}  are devoted to the demonstration, for $x\in A\subset S(p)$ and $i\in \mathbb{N}$, of relation
\begin{equation}\label{convc}
\e(\log\phi_{m,l}(i)|X_i=x) =\e(\log\phi_{m,l}(1)|X_1=x)\to \psi(l) - \log{V_d} -\log q(x), \;\;m\to \infty,
\end{equation}
where $A$ depends on $p$ and $q$ versions, $\p_{X}(S(p)\setminus A)=0$.
Then \textit{Step 4} justifies the desired result \eqref{mr}. \textit{Step~5} contains the validation of \textit{Statement 2}.

{\it Step 1}. Here we establish the distribution convergence for the auxiliary random variables.
Fix any $i\in \mathbb{N}$
and $l\in \{1, \ldots, m\}$.
To simplify notation we do not indicate the dependence of functions on $d$.
For $x\in \mathbb{R}^d$ and $u>0$, we study the asymptotic behavior (as $m\to \infty$) of the following function
\begin{align}\label{eq1a}
\begin{gathered}
    F_{m,l,x}^{i}(u) := \p\left(\phi_{m,l}(i) \leq u | X_i = x\right) = \p \left( m
    V^d_{m,l}(i) \leq u | X_i = x \right) \\ =
    1 - \p \left(V_{m,l}(i) > \left( \frac{u}{m} \right)^{\frac{1}{d}} \Big| X_i = x
    \right) =
    1- \p \left(\norm{x - Y_{(l)}(x,\mathbb{Y}_m)} >\left( \frac{u}{m} \right)^{\frac{1}{d}}\right)\\
    =1 - \sum_{s=0}^{l-1}\binom{m}{s} \left( W_{m,x}(u) \right)^s \left( 1 - W_{m,x}(u) \right)^{m-s} := \p \left(\xi_{m,l,x} \leq u\right),
\end{gathered}
\end{align}
\vspace{-0.2cm}
where
\vspace{-0.2cm}
\begin{equation}\label{p_int}
    W_{m,x}(u) := \int_{B(x, r_m(u))} q(z) \, dz,\;\; r_m(u):=     {\left(\frac{u}{m}\right)}^{\frac{1}{d}},\;\;\xi_{m,l,x} := m \norm{x - Y_{(l)}(x,\mathbb{Y}_m)}^d.
\end{equation}

We have employed in \eqref{eq1a} the independence of random vectors $Y_1,\ldots,Y_m, X_i$ and condition that $Y_1, \ldots, Y_m$ have the same law as $Y$.
We also took into account that an event $\left\{\norm{x - Y_{(l)}(x,\mathbb{Y}_m)} > r_m(u)\right\}$ is
a union of pair-wise disjoint events $A_s$, $s=0,\ldots,l-1$. Here $A_s$ means that
exactly $s$ observations among $\mathbb{Y}_m$ belong to the ball $B(x,r_m(u))$ and other
$m-s$ are outside this ball (probability that $Y$ belongs to the sphere
$\{z\in \mathbb{R}^d:\|z-x\|=r\}$ equals $0$ since $Y$ has a density w.r.t. the Lebesgue measure $\mu$).
Formulas \eqref{eq1a} and \eqref{p_int} show that $F_{m,l,x}^i(u)$ is
the regular conditional distribution function of $\phi_{m,l}(i)$ given $X_i=x$. Moreover, \eqref{eq1a} means that $\phi_{m,l}(i)$, $i \in \{ 1, \ldots, n \}$ are identically distributed and we may omit the dependence on $i$. So, one can replace $F_{m,l,x}^i(u)$ with $F_{m,l,x}(u)$.

According to the Lebesgue differentiation theorem  (see, e.g., \cite{Yeh}, p. 654)    if $q\in L^1(\mathbb{R}^d)$ then, for $\mu$-almost all  $x\in \mathbb{R}^d$,
the following relation holds
\begin{equation}\label{2}
\lim_{r\to 0+}\frac{1}{\mu(B(x,r))}\int_{B(x,r)}|q(z)-q(x)|\,dz=0.
\end{equation}
Let $\Lambda(q)$ stand for a set of all the Lebesgue points of a function $q$,
i.e. points $x\in \mathbb{R}^d$ satisfying \eqref{2}. Clearly, $\Lambda(q)$ depends on the chosen version of $q$ belonging to the class of equivalent functions from $L^1(\mathbb{R}^d)$ and, for an arbitrary version of $q$, we have $\mu(\mathbb{R}^d\setminus \Lambda(q))=0$.

Note that, for each $u>0$, $r_m(u) \to 0$ as $m \to \infty$, and $\mu(B(x,r_m(u))) = V_d {\big(r_m(u)\big)}^d = \frac{V_d u}{m}$. Therefore by virtue of  \eqref{2}, for any fixed $x\in \Lambda(q)$ and $u> 0$,
$$
W_{m,x}(u)  = \frac{V_d\, u}{m} \left(q(x) + \alpha_m{(x,u)}\right),
$$
where $\alpha_m{(x,u)} \to 0, \; m \to \infty$.
Hence, for $x\in \Lambda(q)\cap S(q)$ (thus $q(x)>0$),
due to \eqref{eq1a}
\begin{align}\label{conv}
\begin{gathered}
F_{m,l,x}(u)
\to  1 - \sum_{s=0}^{l-1} \frac{(V_d u q(x))^s}{s!} e^{-V_d u q(x)} := F_{l,x}(u), \;\; m \to \infty.
\end{gathered}
\end{align}
Relation \eqref{conv} means that
\begin{equation}\label{claw}
\xi_{m,l,x}\stackrel{law}\rightarrow \xi_{l,x},\;\;x\in  \Lambda(q)\cap S(q),\;\;m \to \infty,
\end{equation}
where $\xi_{l,x}$ has $\Gamma(V_d\,q(x),l)$ distribution.

We assume without loss of generality (w.l.g.) that, for all $x\in S(q)$, the random variables $\xi_{l,x}$ and $\{\xi_{m,l,x}\}_{m\geq l}$
are defined on a probability space $(\Omega,\mathcal{F},\p)$ since in view of the Lomnicki - Ulam theorem (see, e.g. \cite{Kallenberg}, p. 93) one can consider
the independent copies of $Y_1,Y_2,\ldots$ and $\{\xi_{l,x}\}_{x\in S(q)}$ defined on a certain probability space.
The convergence in law of random variables is preserved under continuous mapping. Hence, for any
$x\in  \Lambda(q)\cap S(q)$, we come to the relation
\begin{equation}\label{b2}
\log \xi_{m,l,x}\stackrel{law}\rightarrow \log \xi_{l,x},\;\;m\to \infty.
\end{equation}
We took into account that, for each $x\in \Lambda(q)\cap S(q)$, one has $\xi_{l,x}>0$ a.s. and
since $Y$ has a density we infer that
$\p(\xi_{m,l,x}>0) =  \p(\norm{x-Y_{(l)}(x, \mathbb{Y}_m)} > 0) = 1$.
More precisely, we can ignore zero values of nonnegative random variables (having zero values with probability zero) when we take their logarithms.

{\it Step 2}. Now we show that instead of \eqref{convc} validity one can verify the following statement. For $\mu$-almost every $x\in \Lambda(q) \cap  S(q)$,
\begin{equation}\label{a2}
\e \log \xi_{m,l,x}\to \e \log \xi_{l,x},\;\;m\to \infty.
\end{equation}
Note that if $\eta\sim \Gamma(\alpha, \lambda)$, where $\alpha > 0$ and $\lambda >0$, then
\begin{align}\label{calculate_integral_erlang}
\begin{gathered}
    \e \log \eta = \int_{(0,\infty)} \log{u} \, \frac{\alpha^\lambda u^{\lambda-1} e^{-\alpha u}}{\Gamma(\lambda)} \,du =
    \int_{(0, \infty)} \left( \log{\frac{v}{\alpha}} \right) \frac{v^{\lambda-1} e^{-v}}{\Gamma{(\lambda)}} \, dv\\ = \int_{(0, \infty)} \log{v} \frac{v^{\lambda-1} e^{-v}}{\Gamma{(\lambda)}} \, dv - \log{\alpha} \frac{\int_{(0, \infty)} v^{\lambda-1} e^{-v} \, dv}{\Gamma(\lambda)} =  \psi(\lambda) - \log{\alpha}.
\end{gathered}
\end{align}
Set $\alpha = V_d q(x)$, where $q(x)>0$ for $x\in S(q)$, and $\lambda = l$. Then
$\e \log \xi_{l,x}= \psi{(l)} - \log{(V_d q(x))}= \psi{(l)} - \log{V_d} - \log{q(x)}$.
By virtue of \eqref{ce1}, for each $x\in \mathbb{R}^d$,
\begin{gather*}
    \e \log{\xi_{m,l,x}} = \int_{(0,\infty)} \log{u} \, dF_{m,l,x}(u) = \int_{(0,\infty)} \log{u} \, d
    \p(\phi_{m,l}(1) \leq u|X_1 = x)  \\
    = \e(\log{\phi_{m,l}(1)} | X_1 = x).
\end{gather*}
Thus, for $x\in \Lambda(q) \cap S(q)$, the relation $\e(\log\phi_{m,l}(1)) | X_1 = x) \to \psi{(l)} - \log{V_d} - \log{q(x)}$ holds if and only if \eqref{a2} is true.

According to Theorem 3.5 \cite{Billingsley} we would have established \eqref{a2} if relation \eqref{b2} could be
supplemented, for $\mu$-almost all $x\in \Lambda(q) \cap  S(q)$, by the uniform integrability of a family $\{\log \xi_{m,l,x}\}_{m\geq m_0(x)}$.
Note that, for each $N\in \mathbb{N}$, a function $G_N(t)$ introduced by \eqref{GN} is increasing on $(0,\infty)$ and $\frac{G_N(t)}{t}\to \infty$, as $t\to \infty$. Therefore,
by the de la Valle Poussin theorem (see, e.g., Theorem 1.3.4 \cite{Borkar}),
to guarantee, for $\mu$-almost every $x\in \Lambda(q) \cap  S(q)$, the uniform integrability of $\{\log \xi_{m,l,x}\}_{m\geq m_0(x)}$
it suffices to prove, for such $x$, a positive $C_0(x)$ and $m_0(x)\in \mathbb{N}$, that
\begin{equation}\label{b3}
\sup_{m\geq m_0(x)} \e G_{N_1}(|\log \xi_{m,l,x}|)\leq C_0(x)<\infty,
\end{equation}
where $G_{N_1}$ appears in conditions of Theorem \ref{th1}.

{\it Step 3} is devoted to proving validity of \eqref{b3}.
It is convenient to divide this proof into its own parts (3a), (3b), etc.
For any  $N \in \mathbb{N}$, set
$$
g_N(t)=
\begin{cases}
-\frac{1}{t}\left(\log_{[N]}(-\log t)+ \frac{1}{\prod_{j=1}^{N-1}
\log_{[j]}(-\log t)}\right), & t \in \left(0,\frac{1}{e_{[N]}}\right],\\
0, & t \in  \left(\frac{1}{e_{[N]}},e_{[N]}\right],\\
\frac{1}{t}\left(\log_{[N]}(\log{t}) +  \frac{1}{\prod_{j=1}^{N-1} \log_{[j]}(\log{t})}\right),
& t \in \left(e_{[N]}, \infty\right),
\end{cases}
$$
where the product over empty set (when $N=1$) is equal to 1.

We will employ the following result,
its proof is given in Appendix.

\begin{lem}\label{lemma_G}
Let $F(u), u\in \mathbb{R}$, be a distribution function such that $F(0)=0$. Then, for each $N\in \mathbb{N}$, one has

1) $\int_{\left(0,\frac{1}{e_{[N]}}\right]}G_N(|\log u|) dF(u) = \int_{\left(0,\frac{1}{e_{[N]}}\right]}F(u)(-g_N(u)) du$,

2) $\int_{\left(e_{[N]}, \infty\right)}G_N(|\log u|) dF(u) = \int_{\left(e_{[N]}, \infty\right)}(1-F(u)) g_N(u) du$.
\end{lem}

\noindent
Note that, for $u \in \left(\frac{1}{e_{[N_1]}},e_{[N_1]}\right]$, we have
$G_{N_1}(|\log u|)=0$.
Therefore, due to Lemma~\ref{lemma_G}, for $x\in \Lambda(q) \cap  S(q)$ and $m \geq l$, we get
$\e G_{N_1}(|\log{\xi_{m,l,x}}|):=I_1(m,x)+I_2(m,x)$ where
\begin{gather*}
    I_1(m,x):= \int_{\left(0,\frac{1}{e_{[N_1]}}\right]}F_{m,l,x}(u)(-g_{N_1}(u))du,\;\;\;
    I_2(m,x):=\int_{(e_{[N_1]},\infty)}(1-F_{m,l,x}(u))g_{N_1}(u)du.
    \end{gather*}
For convenience sake we write $I_1(m,x)$ and $I_2(m,x)$ without indicating their dependence on $N_1, l$ and $d$. Recall that  $N_1$ is fixed.

{\it Part (3a)}. We provide bounds for $I_1(m,x)$.
Take  $R_1>0$ appearing in conditions of  Theorem \ref{th1} and any $u\in \left(0,\frac{1}{e_{[N_1]}}\right]$. Let us denote $m_1 := \max\left\{ \left\lceil\frac{1}{e_{[N_1]} R_1^d}\right\rceil, l \right\}$, where $\lceil a \rceil := \inf\{m\in \mathbb{Z}: m\geq a\}$, $a\in \mathbb{R}$. Then
$
r_m(u) = \left( \frac{u}{m} \right)^{1/d} \leq {\left( \frac{1}{e_{[N_1]} m} \right)}^{1/d} \leq R_1
$
if $m\geq m_1$. Note also that we can consider only $m \geq l$ everywhere below, because the size of sample $\mathbb{Y}_m$ should not be less than number of the neighbors $l$ (see, e.g., \eqref{eq1a}). Thus, for $R_1>0$, $u\in \left(0,\frac{1}{e_{[N_1]}}\right]$, $x\in \mathbb{R}^d$ and $m\geq m_1$,
\begin{gather*}
    \frac{W_{m,x}(u)}{\mu(B(x, r_m(u)))} = \frac{\int_{B(x,r_m(u))} q(y) \, dy}{r_m^d(u) V_d} \leq \sup_{r \in (0, R_1]} \frac{\int_{B(x,r)} q(y) \, dy}{r^d V_d} = M_q(x,R_1),
\end{gather*}
and we obtain an inequality
\begin{equation}\label{third: eq1}
W_{m,x}(u) \leq M_q(x, R_1) \, \mu(B(x, r_m(u))) =  \frac{M_q(x, R_1) V_d\, u }{m}.
\end{equation}
If $\varepsilon \in (0,1]$ and $t \in [0,1]$ then, for all $m \geq 1$,
invoking the Bernoulli inequality, one has
\begin{equation}\label{Bernoulli}
1-(1-t)^m \leq (mt)^{\varepsilon}.
\end{equation}
By assumptions of the Theorem  $Q_{p,q}(\varepsilon_1, R_1) < \infty$ for some $\varepsilon_1 > 0$, $R_1>0$. According to Lemma \ref{lemma1} we can assume that $\varepsilon_1 < 1$. Thus, due to  \eqref{Bernoulli} and since $W_{m,x}(u) \in [0,1]$ for all $x \in \mathbb{R}^d$, $u > 0$ and $m \geq l$, we get
\begin{equation}\label{apply_lemma3}
1-(1-W_{m,x}(u))^m \leq (m W_{m,x}(u))^{\varepsilon_1}.
\end{equation}
In view of \eqref{eq1a}, \eqref{third: eq1} and \eqref{apply_lemma3} one can claim now that, for all $x \in \Lambda(q) \cap S(q)$, $u \in (0,\frac{1}{e_{[N]}}]$  and $m \geq m_1$,
\begin{align}\label{ineq_B}
\begin{gathered}
    F_{m,l,x}(u) = 1 - \sum_{s=0}^{l-1}\binom{m}{s} \left( W_{m,x}(u) \right)^s \left( 1 - W_{m,x}(u) \right)^{m-s}  \\ \leq 1-(1-W_{m,x}(u))^{m} \leq \left(m \frac{M_{q}(x, R_1) V_d u}{m}\right)^{\varepsilon_1} = (M_{q}(x, R_1))^{\varepsilon_1} V_d^{\varepsilon_1} u^{\varepsilon_1}.
\end{gathered}
\end{align}
Therefore, for any $x\in \Lambda(q) \cap S(q)$ and $m\geq m_1$, one can write
\begin{align}\label{third: eq2}
\begin{gathered}
I_1(m,x) \leq (M_q(x, R_1))^{\varepsilon_1} V_d^{\varepsilon_1} \int_{\left(0,\frac{1}{e_{[N_1]}}\right]} u^{\varepsilon_1}(-g_{N_1}(u))  \, du \\
\leq (M_q(x, R_1))^{\varepsilon_1} V_d^{\varepsilon_1} \int_{\left(0,\frac{1}{e_{[N_1]}}\right]} \frac{\log_{[N_1]}(-\log u) +1}{u^{1-\varepsilon_1}} du = U_1(\varepsilon, N, d) (M_q(x, R_1))^{\varepsilon_1},
\end{gathered}
\end{align}
where $U_1(\varepsilon, N, d) := V_d^{\varepsilon} L_N(\varepsilon)$, $L_N(\varepsilon) := \int_{[e_{[N-1]}, \infty)} (\log_{[N]}(t) +1) e^{-\varepsilon t} dt < \infty$ for each $\varepsilon > 0$ and any $N \in \mathbb{N}$. We took into account
that $(-g_{N_1}(u))\leq \frac{1}{u}(\log_{[N_1]}(-\log u) +1)$ if $u\in \left(0,\frac{1}{e_{[N_1]}}\right]$.

{\it Part (3b)}. We give bounds for $I_2(m,x)$. Since $g_{N_1}(u) \leq \frac{\log_{[N_1+1]}(u) + 1}{u}$ if  $u \in (e_{[N_1]}, \infty)$, we can write, for $m \geq \max\{e^2_{[N_1]}, l\}$,
\begin{align*}
\begin{gathered}
I_2(m,x) \leq
    \int_{(e_{[N_1]}, {\sqrt{m}}]} (1 - F_{m,l,x}(u)) \frac{\log_{[N_1+1]}(u) + 1}{u} \, du \,\, \\+ \int_{({\sqrt{m}},{m^2}]} (1 - F_{m,l,x}(u)) \frac{\log_{[N_1+1]}(u) + 1}{u} \, du + \int_{\left(m^2, \, \infty \right]} (1-F_{m,l,x}(u)) g_{N_1}(u) \, du \\ := J_1(m,x) + J_2(m,x) + J_3(m,x).
\end{gathered}
\end{align*}
Evidently,
\begin{align}\label{F_transform}
\begin{gathered}
    1-F_{m,l,x}(u) =
    \sum_{r=m-l+1}^{m} \binom{m}{r} \left( P_{m,x}(u) \right)^r \left( 1-P_{m,x}(u) \right)^{m-r} = \p (Z \geq m-l+1),
\end{gathered}
\end{align}
where $P_{m,x}(u) = 1-W_{m,x}(u)$ and $Z \sim {\sf Bin}(m,P_{m,x}(u))$.

By Markov's inequality $\p(Z \geq x) \leq e^{-\lambda x}\e e^{\lambda Z}$
for any $\lambda > 0$ and $x > 0$.
One has
\begin{align}\label{F_transform22}
\begin{gathered}
    \e e^{\lambda Z} = \sum_{j=0}^m e^{\lambda j} \binom{m}{j} \left( P_{m,x}(u) \right)^j \left( 1-P_{m,x}(u) \right)^{m-j}  \\ =\sum_{j=0}^m \binom{m}{j} \left( P_{m,x}(u) e^{\lambda} \right)^j \left( 1-P_{m,x}(u) \right)^{m-j} = \big(1-P_{m,x}(u) + e^{\lambda} P_{m,x}(u)\big)^m.
\end{gathered}
\end{align}
Consequently, for each $\lambda >0$,
\begin{align}\label{proba_calculation}
\begin{gathered}
    1-F_{m,l,x}(u) \leq e^{-\lambda(m-l+1)}\big(1-P_{m,x}(u) + e^{\lambda} P_{m,x}(u)\big)^m \\= e^{-\lambda(m-l+1)}\big(W_{m,x}(u) + e^{\lambda} (1-W_{m,x}(u))\big)^m  =
    e^{\lambda (l-1)} \left( 1 - \left( 1-\frac{1}{e^{\lambda}} \right) W_{m,x}(u) \right)^m.
\end{gathered}
\end{align}
To simplify bounds we take $\lambda = 1$ and set $S_1 = S_1(l):= e^{l-1}$, $S_2 := 1- \frac{1}{e}$
(recall that $l$ is fixed). Thus $S_1 \geq 1$ and $S_2 < 1$. Therefore,
\begin{equation}\label{main_1-f_calculation}
    1-F_{m,l,x}(u) \leq S_1 \left( 1 - S_2 \, W_{m,x}(u) \right)^m \leq
    S_1 \exp\left\{  - S_2 \, m  W_{m,x}(u) \right\},
\end{equation}
where we have used an elementary inequality $1-t \leq e^{-t}$, $t\in [0,1]$.

For $R_2 > 0$ appearing in conditions of the Theorem and any $u\in \left(e_{[N]},\sqrt{m}\right]$, one can choose $m_2 := \max\left\{\left\lceil\frac{1}{R_2^{2d}}\right\rceil, \left\lceil e_{[N_1]}^2 \right\rceil, l\right\}$ such that if $m\geq m_2$ then
$
r_m(u) = \left( \frac{u}{m} \right)^{1/d} \leq \left(\frac{1}{\sqrt{m}}\right)^{1/d} \leq R_2.
$
Due to \eqref{p_int} and  \eqref{main_1-f_calculation}, for $u \in (e_{[N_1]}, \sqrt{m}]$ and $m \geq m_2$, one has
\begin{align}\label{1-F_add}
\begin{gathered}
        1 - F_{m,l,x}(u) \leq S_1 \exp\left\{ -S_2 \, m \frac{V_d u}{m} \frac{W_{m,x}(u)}{\frac{V_d u}{m}} \right\}  \\=
        S_1 \exp\left\{ -S_2 V_d u \frac{\int_{B(x,r_m(u))} q(z) \, dz}{\mu(B(x, r_m(u)))} \right\} \leq S_1 \exp\left\{ -S_2 V_d u \, m_q(x, R_2)\right\},
\end{gathered}
\end{align}
by definition of $m_f$ (for $f=q$) in \eqref{mM}.
Now we use the following Lemma 3.2 of \cite{Bul_Dim}.
\begin{lem}\label{l1}
For a version of a density $q$ and each $R>0$, one has $\mu(S(q)\setminus D_q(R))=0$ where $D_q(R):=\{x\in S(q): m_q(x,R)>0\}$ and $m_q(\cdot,R)$ is defined according to \eqref{mM}.
\end{lem}

It is easily seen that, for any $t>0$ and each $\delta \in (0,e]$, one has $e^{-t}\leq t^{-\delta}$.
Thus, for $x\in D_q(R_2)$, $m\geq m_2$,  $u\in (e_{[N]},\sqrt{m}]$ and $\varepsilon_2 > 0$, we deduce from conditions of the Theorem (in view of Lemma \ref{lemma1} one can suppose that $\varepsilon_2 \in (0, e]$),
taking into account that $m_q(x,R_2)>0$ for $x\in D_q(R_2)$ and applying  relation \eqref{1-F_add}, that
\begin{equation}\label{1-F}
1 - F_{m,l,x}(u)\leq S_1 \left( S_2 V_d u \, m_q(x, R_2) \right)^{-\varepsilon_2}.
\end{equation}
Thus, for all $x \in \Lambda(q)\cap S(q) \cap D_q(R_2)$ and any $m \geq m_2$,
\begin{align}\label{estJ1}
\begin{gathered}
    J_1(m,x)
     \leq \frac{S_1}{(S_2 \, V_d)^{\varepsilon_2}(m_q(x, R_2))^{\varepsilon_2}}
     \int_{(e_{[N_1]}, \infty)}\!\!\!
     \frac{\log_{[N_1+1]}(u) + 1}{u^{1+\varepsilon_2}} \, du\\
          = U_2(\varepsilon_2, N_1, d, l) (m_q(x, R_2))^{-\varepsilon_2},
\end{gathered}
\end{align}
where $U_2(\varepsilon, N, d, l) := S_1(l) \, L_N(\varepsilon)(S_2 \, V_d)^{-\varepsilon}$.

{\it Part (3c)}. Consider $J_2(m,x)$.
In view of \eqref{1-F}, for all $x \in \Lambda(q)\cap S(q) \cap D_q(R_2)$ and any $m \geq m_2$, it holds $1 - F_{m,l,x}(\sqrt{m})\leq S_1 \left( S_2 V_d \, m_q(x, R_2) \sqrt{m} \right)^{-\varepsilon_2}$. Thus (as $m_2 \geq 2$)
\begin{align*}
\begin{gathered}
    J_2(m,x) \leq \int_{\left(\sqrt{m}, \,m^2\right]} (1-F_{m,l,x}(u)) \frac{\log_{[N_1+1]}(u) + 1}{u} \, du  \\ \leq\left(1-F_{m,l,x}(\sqrt{m})\right) \int_{\left(\sqrt{m}, \, m^2\right]} \left(\log_{[N_1+1]}(u) + 1\right)\, d \log{u}  \\
    \leq S_1 (S_2 V_d)^{-\varepsilon_2} \, \left(m_q(x, R_2)\right)^{-\varepsilon_2} m^{-\frac{\varepsilon_2}{2}} \left( \log_{[N_1]}(2 \log{m}) + 1\right) \frac{3}{2}\log{m}.
\end{gathered}
\end{align*}
Then, for all $x \in \Lambda(q)\cap S(q) \cap D_q(R_2)$ and any $m \geq m_2$,
\begin{equation}\label{U3_inequality}
    J_2(m,x) \leq U_3(m, \varepsilon_2, N_1, d, l) \left(m_q(x, R_2)\right)^{-\varepsilon_2},
\end{equation}
where $U_3(m, \varepsilon_2, N_1, d, l):= \frac{3}{2} S_1(l) (S_2 V_d)^{-\varepsilon_2} m^{-\frac{\varepsilon_2}{2}} \log{m} \left( \log_{[N_1]}(2 \log{m}) + 1 \right) \to 0$, $m \to \infty$.

\vspace{0.2cm}
{\it Part (3d)}. To get bounds for $J_3(m,x)$ we employ several auxiliary results.

\begin{lem}\label{l4}
For each $N \in \mathbb{N}$ and any $\nu > 0$, there are  $a := a(d, \nu) \geq 0, \, b := b(N,d, \nu)\geq 0$ such that, for arbitrary $x,y\in \mathbb{R}^d$,
\begin{equation}\label{v0}
G_{N}\left(|\log{\|x-y\|^d}|^{\nu}\right) \leq a \, G_N\left(|\log{\|x-y\|}|^{\nu}\right) + b.
\end{equation}
\end{lem}

The proof is provided in Appendix.

On the one hand, by \eqref{p_int}, for any $w \geq 0$, we get
$$
W_{m,x}(m w) = \int_{B(x, w^{1/d})} q(z) \, dz = W_{1,x}(w).
$$
On the other hand, by \eqref{eq1a}, one has $F_{1,1,x}(w) = 1 - \big(1-W_{1,x}(w)\big) = W_{1,x}(w)$. Consequently, for any $m \in \mathbb{N}$, $w \geq 0$ and all $x \in \mathbb{R}^d$,
\begin{equation}\label{obs_2}
W_{m,x}(m w) = F_{1,1,x}(w).
\end{equation}
Moreover, $F_{1,1,x}(w) = \p(\norm{Y-x}^d \leq w)$. So, $\xi_{1,1,x}\stackrel{law}=\norm{Y-x}^d$. Thus, in view of Lemmas \ref{lemma_G} and \ref{l4} (for $N=N_1$ and $\nu = 1$)
\begin{align}\label{link_xi_to_condition_one}
\begin{gathered}
    \int_{\left(e_{[N_1]}, \infty \right)} (1-F_{1,1,x}(w)) g_{N_1}(w) \, dw = \int_{\left(e_{[N_1]}, \infty \right)} G_{N_1}(\log{w}) \, d F_{1,1,x}(w)  \\=
    \e \left[G_{N_1}\left(\log \xi_{1,1,x}\right) \ind\left\{\xi_{1,1,x} > e_{[N_1]}\right\} \right] = \e [G_{N_1}(\log \|Y-x\|^d) \ind\{\|Y-x\|^d > e_{[N_1]}\} ]  \\=
    \int_{y\in \mathbb{R}^d, \,\norm{x-y} > \left(e_{[N_1]}\right)^{1/d}} G_{N_1}(\log \norm{x-y}^d) q(y) \, dy  \\\leq
        a(d,1) \int_{y\in \mathbb{R}^d, \,\norm{x-y} > \left(e_{[N_1]}\right)^{1/d}} G_{N_1}(|\log \norm{x-y}|) q(y) \, dy + b(N_1,d,1)  \\=
        a(d,1) \int_{y\in \mathbb{R}^d, \,\norm{x-y} > e_{[N_1]}} G_{N_1}(\log \norm{x-y}) q(y) \, dy + b(N_1,d,1),
\end{gathered}
\end{align}
since $G_{N}(t)=0$ for $t\in [0,e_{[N-1]}]$, $N\in \mathbb{N}$.

Now we will estimate $1 - F_{m,l,x}(u)$ in a way  different from
\eqref{proba_calculation}.
Fix any $\delta > 0$. Note that, for all $m \geq (l-1) \left( 1 + \frac{1}{\delta} \right)$ and $s \in \{0, \ldots, l-1\}$, it holds $\frac{m}{m-s} \leq \frac{m}{m-l+1} \leq 1+\delta$. Then, for all $x \in \mathbb{R}^d$, $u \geq 0$ and $m \geq \max\{l, (l-1) \left( 1 + \frac{1}{\delta} \right)\}$, in view of \eqref{eq1a} one can write
\begin{align*}
\begin{gathered}
    1 - F_{m,l,x}(u)
        =  \left( 1 - W_{m,x}(u) \right) \sum_{s=0}^{l-1}\binom{m-1}{s} \frac{m}{m-s} \left( W_{m,x}(u) \right)^s \left( 1 - W_{m,x}(u) \right)^{(m-1)-s}  \\
    \leq (1 + \delta) \left( 1 - W_{m,x}(u) \right) \sum_{s=0}^{l-1}\binom{m-1}{s} \left( W_{m,x}(u) \right)^s \left( 1 - W_{m,x}(u) \right)^{(m-1)-s}
\end{gathered}
\end{align*}
\vspace{-0.2cm}
\begin{equation}\label{way2_not_full}
      \leq (1 + \delta) \left( 1 - W_{m,x}(u) \right).
\end{equation}

We are going to employ the following statement as well.
\begin{lem}\label{sv}
For each $N\in \mathbb{N}$, a function $\log_{[N]}(t)$, $t > e_{[N-1]}$,
is slowly varying at infinity.
\end{lem}

Its proof is elementary and thus is omitted.

{\it Part (3e)}. Now we are ready to get the bound for $J_3(m,x)$. Set $u = m w$. Then one has
\begin{align*}
    \begin{gathered}
      J_3(m,x) = \int_{\left(m^2, \, \infty \right]} (1-F_{m,l,x}(u)) \frac{1}{u}\left(\log_{[N_1]}(\log{u}) +  \frac{1}{\prod_{j=1}^{N_1-1} \log_{[j]}(\log{u})}\right) \, du  \\
      =\int_{\left(m, \, \infty \right]} (1-F_{m,l,x}(mw)) \frac{1}{w}\left(\log_{[N_1+1]}(m w) +  \frac{1}{\prod_{j=2}^{N_1} \log_{[j]}(m w)}\right) \, dw.
    \end{gathered}
\end{align*}
Inequality $w > m$ and Lemma \ref{sv} imply $\log_{[N_1+1]}(m w) \leq \log_{[N_1+1]}(w^2) = \log_{[N_1]}(2 \log{w}) \leq 2 \log_{[N_1+1]}(w)$ for $w$ large enough, namely for all $w \geq W$, where $W=W(N_1)$.

Take $\delta > 0$ and set $m_3 := \max\left\{ l, \left\lceil (l-1) \left( 1 + \frac{1}{\delta}\right) \right\rceil, \left\lceil W(N_1) \right\rceil, \left\lceil e_{[N_1]} \right\rceil \right\}$. Let further $m \geq m_3$. Then
\begin{align*}
    \begin{gathered}
      J_3(m,x) \leq 2 \int_{\left(m, \, \infty \right]} (1-F_{m,l,x}(mw)) \frac{1}{w}\left( \log_{[N_1+1]}(w) +  \frac{1}{\prod_{j=2}^{N_1} \log_{[j]}(w)}\right) \, dw.
    \end{gathered}
\end{align*}
By virtue of \eqref{obs_2} and \eqref{way2_not_full}  one has
\begin{equation}\label{another_way_to_estimate_with_mw}
1-F_{m,l,x}(mw) \leq (1+\delta) \left( 1-W_{m,x}(mw) \right) = (1+\delta) \left( 1 - F_{1,1,x}(w) \right).
\end{equation}
Hence it can be seen that
\begin{align}\label{J_3_last_in_theorem1}
    \begin{gathered}
        J_3(m,x)
        \leq 2 (1+\delta) \int_{\left(m, \, \infty \right]} (1-F_{1,1,x}(w)) g_{N_1}(w) \, dw.
    \end{gathered}
\end{align}
Introduce
$$
R_{N}(x) := \int_{y\in \mathbb{R}^d,\,\norm{x-y} > e_{[N]}} G_{N}(\log{\norm{x-y}}) q(y) \, dy,\;\;A_p(G_N):=\{x\in S(p): R_{N}(x) < \infty\}.
$$

Let us note: 1) $\p_{X}(S(p)\setminus A_p(G_{N_1}))=0$ as we assumed that $K_{p,q}(1,N_1)<\infty$;
2) $\p_{X}(S(p) \setminus S(q)) = 0$ as $\p_X\ll \p_Y$;
3) $\mu\big(S(q) \setminus (\Lambda(q) \cap D_q(R_2))\big) = 0$ due
to Lemma \ref{l1}.
Since $\p_X \ll \mu$ we conclude that $\p_X \big(S(q) \setminus (\Lambda(q) \cap D_q(R_2)) \big) = 0$. Hence, one has $\p_X\big(S(p) \setminus (\Lambda(q) \cap D_q(R_2))\big) = 0$ in view of 2) and because $B \setminus C \subset (B \setminus A) \cup (A \setminus C)$  for any $A,B,C \subset \mathbb{R}^d$. Set further $A:=\Lambda(q) \cap S(q) \cap D_q(R_2)\cap S(p) \cap A_p(G_{N_1})$.
It follows from 1), 2) and 3) that $\p_X(S(p) \setminus A)=0$, so $\p_X(A) = 1$. We are going to consider only $x \in A$.

Then, by virtue of \eqref{link_xi_to_condition_one} and \eqref{J_3_last_in_theorem1}, for all $m \geq m_3$ and $x \in A$, we come to the inequality
\begin{equation}\label{J22_mx}
    J_3(m,x) \leq 2 (1 + \delta) \big(a(d,1) R_{N_1}(x) + b(N_1, d,1)\big)
    = A(\delta, d) R_{N_1}(x) + B(\delta, d, N_1),
\end{equation}
where $A(\delta, d) := 2 (1+\delta) a(d,1)$, $B(\delta, d, N_1) := 2 (1 + \delta) b(N_1, d, 1)$.

{\it Part (3f)}. Thus, for each $x\in A$ and $m \geq \max\{ m_1, m_2, m_3 \}$, taking into account \eqref{third: eq2}, \eqref{estJ1}, \eqref{U3_inequality} and \eqref{J22_mx} we can claim that
\begin{align}\label{final_inequality_!}
    \begin{gathered}
        \e G_{N_1}(|\log{\xi_{m,l,x}}|) \leq I_1(m,x) + J_1(m,x) + J_2(m,x) + J_3(m,x) \\ \leq U_1(\varepsilon_1, N_1, d) (M_q(x, R_1))^{\varepsilon_1} + U_2(\varepsilon_2, N_1, d, l) (m_q(x, R_2))^{-\varepsilon_2}  \\ + \,U_3(m, \varepsilon_2, N_1, d, l) \left(m_q(x, R_2)\right)^{-\varepsilon_2} + \left(A(\delta, d) R_{N_1}(x) + B(\delta, d, N_1)\right).
    \end{gathered}
\end{align}
Moreover, for any $\kappa > 0$, one can take $m_4 = m_4(\kappa, \varepsilon_2, N_1, d, l)\in \mathbb{N}$ such that $U_3(m, \varepsilon_2, N_1, d, l) \leq \kappa$ for $m \geq m_4$. Then by virtue of \eqref{final_inequality_!}, for each $x\in A$ and $m \geq m_0 := \max\{ m_1, m_2, m_3, m_4 \}$,
\begin{align}\label{final_inequality_kappa}
    \begin{gathered}
        \e G_{N_1}(|\log{\xi_{m,l,x}}|) \leq  U_1(\varepsilon_1, N_1, d) (M_q(x, R_1))^{\varepsilon_1}  \\ +\big(U_2(\varepsilon_2, N_1, d, l) + \kappa\big) (m_q(x, R_2))^{-\varepsilon_2}  + \left(A(\delta, d) R_{N_1}(x) + B(\delta, d, N_1)\right) := C_0(x) < \infty.
    \end{gathered}
\end{align}
Hence, for each $x \in A$, the uniform integrability of the family $\left\{\log{\xi_{m,l,x}}\right\}_{m \geq m_0}$ is established.

{\it Step 4}. Now we verify \eqref{main1}. We have already proved, for each $x\in A$ (thus, for $\p_{X}$-almost every $x$ belonging to  $S(p)$) that $\e(\log\phi_{m,l}(1)|X_1=x) \to \psi(l) - \log{V_d} -\log q(x)$, $m\to \infty$.
Set $Z_{m,l}(x):=\e(\log\phi_{m,l}(1)|X_1=x) = \e \log{\xi_{m,l,x}}$.
Consider $x\in A$ and take any $m\geq \max\{m_1,m_2, m_3,m_4\}$. We use the following property of $G_N$ which is shown in Appendix.
\begin{lem}\label{l6}
For each $N\in \mathbb{N}$, a function $G_N$ is convex on $\mathbb{R}_+$.
\end{lem}

Thus a function $G_{N_1}$ is nondecreasing and convex.
On account of the Jensen inequality
\begin{equation}\label{fourth: eq1}
\begin{gathered}
    G_{N_1}(|Z_{m,l}(x)|) = G_{N_1}(|\e \log{\xi_{m,l,x}}|) \leq G_{N_1}(\e |\log{\xi_{m,l,x}}|)
     \leq  \e G_{N_1}(|\log{\xi_{m,l,x}}|).
\end{gathered}
\end{equation}
Relation \eqref{final_inequality_kappa} guarantees that, for all $m \geq m_0$,
\begin{gather*}
    \int_{\mathbb{R}^d} G_{N_1}(|Z_{m,l}(x)|) p(x) \, dx \leq
    U_1(\varepsilon_1, N_1, d) Q_{p,q}(\varepsilon_1, R_1)  \\ + \big(U_2(\varepsilon_2, N_1, d, l) + \kappa\big) T_{p,q}(\varepsilon_2, R_2) +  A(\delta, d) K_{p,q}(1, N_1) + B(\delta, d, N_1).
\end{gather*}
We have established  uniform integrability of  the family $\{Z_{m,l}\}_{m \geq m_0}$
w.r.t. measure $\p_{X}$.
Thus, for $i\in \mathbb{N}$,
\begin{gather*}
\e\log\phi_{m,l}(i)= \int_{\mathbb{R}^d}\e(\log\phi_{m,l}(1)|X_1=x)\p_{X_1}(dx) =
\int_{\mathbb{R}^d}Z_{m,l}(x)\;p(x)dx\\
\to \psi(l) - \log{V_d} -\int_{\mathbb{R}^d}p(x)\log q(x)dx,\;\;m\to \infty,
\end{gather*}
and we
come to relation \eqref{mr}.

{\it Step 5}. Let us briefly discuss  the \textit{Statement 2}.
Similar to $F_{m,l,x}(u)$, one can introduce, for $n,k\in \mathbb{N}$, $n\geq k+1$, $x\in \mathbb{R}^d$ and $u>0$, the following function
\begin{align}\label{eq2a}
\begin{gathered}
	\widetilde{F}_{n,k,x}(u) := \p\left(\zeta_{n,k}(i) \leq u | X_i=x\right) = 1 - \p \left( \norm{x- X_{(k)}(x, \mathbb{X}_n \setminus \{x\})} > r_{n-1}(u) \right) \\
   = 1 - \sum_{s=0}^{k-1}\binom{n-1}{s} \left( V_{n-1,x}(u) \right)^s \left( 1 - V_{n-1,x}(u) \right)^{n-1-s} := \p \left(\widetilde{\xi}_{n,k,x} \leq u\right),
\end{gathered}
\end{align}
where $r_n(u)$ was defined  in \eqref{p_int},
\vspace{-0.2cm}
\begin{equation}\label{p_int_2}
    V_{n,x}(u) := \int_{B(x, r_n(u))} p(z) \, dz,\;\;\widetilde{\xi}_{n,k,x} := (n-1) \norm{x - X_{(k)}(x,\mathbb{X}_n \setminus \{x\})}^d.
\end{equation}
Formulas \eqref{eq2a} and \eqref{p_int_2} show that $\widetilde{F}_{n,k,x}(u)$ is
the regular conditional distribution function of $\zeta_{n,k}(i)$ given $X_i=x$. Moreover, for any fixed $u > 0$ and $x \in \Lambda(p) \cap S(p)$ (thus $p(x) > 0$),
\vspace{-0.2cm}
\begin{align*}
\begin{gathered}
\widetilde{F}_{n,k,x}(u)
\to  1 - \sum_{s=0}^{k-1} \frac{(V_d u\, p(x))^s}{s!} e^{-V_d u\, p(x)} := \widetilde{F}_{k,x}(u), \;\; n \to \infty.
\end{gathered}
\end{align*}
Hence, $\widetilde{\xi}_{n,k,x}\stackrel{law}\rightarrow \widetilde{\xi}_{k,x}$, $x\in  \Lambda(p)\cap S(p)$, $n \to \infty$. For $N\in \mathbb{N}$, set $\widetilde{A}_p(G_N):=\{x\in S(p): \widetilde{R}_{N}(x) < \infty\}$, where
$$
\widetilde{R}_{N}(x) := \int_{y\in \mathbb{R}^d,\,\norm{x-y} > e_{[N]}} G_{N}(\log{\norm{x-y}}) p(y)dy.
$$
Introduce $\widetilde{A} := \Lambda(p) \cap S(p) \cap D_p(R_4) \cap \widetilde{A}_p(G_{N_2})$. Then $\p(\widetilde{A})=1$ and, for $x \in \widetilde{A}$, one can verify that
$\e G_{N_2}(|\log{\widetilde{\xi}_{n,k,x}}|) \leq \widetilde{C}_0(x) < \infty$ and therefore
$\e \log \widetilde{\xi}_{n,k,x} \to \e \log \widetilde{\xi}_{k,x}$. Thus $\e(\log\zeta_{n,k}(1)|X_1=x) \to \psi(k) - \log{V_d} -\log p(x)$, $n\to \infty$.
Set $\widetilde{Z}_{n,k}(x) := \e( \log\zeta_{n,k}(1) | X_1 = x )$.
One can see that, for all $n\geq n_0$, $\int_{\mathbb{R}^d} G_{N_2}(|\widetilde{Z}_{n,k}(x)|) p(x) \, dx <\infty$.
Hence similar to Steps 1--4 we come to relation \eqref{mra}.

The proof of Theorem \ref{th1} is complete. $\square$

\section{Proof of Theorem \ref{th_main2}}

First of all note that, in view of
Lemma \ref{lemma1}, the finiteness of $K_{p,q}(2, N_1)$ and $K_{p,p}(2, N_2)$ implies the finiteness of $K_{p,q}(1, N_1)$ and $K_{p,p}(1, N_2)$, respectively. Thus the conditions of Theorem \ref{th_main2} entail validity of Theorem \ref{th1} statements. Consequently under the conditions of Theorem  \ref{th_main2}, for $n$ and $m$ large enough, one can claim that $\widehat{D}_{n,m}(k,l) \in L^1(\Omega)$ and $\e \widehat{D}_{n,m}(k,l) \to D(\p_X || \p_Y)$, as $n,m \to \infty$.

We will show that  $\widehat{D}_{n,m}(k,l) \in L^2(\Omega)$ for all $n$ and $m$ large enough. Then we can write
$$
\e \left(\widehat{D}_{n,m}(k,l) - D(\p_X || \p_Y)\right)^2 = \var\left(\widehat{D}_{n,m}(k,l)\right) + \left(\e \widehat{D}_{n,m}(k,l) - D(\p_X || \p_Y)\right)^2.
$$
Therefore to prove \eqref{main2} we will demonstrate that
$\var \left(\widehat{D}_{n,m}(k,l)\right) \to 0$, $n,m\to \infty$.

Due to \eqref{eq1a} the random variables
$\log{\phi_{m,l}(1)}, \ldots, \log{\phi_{m,l}(n)}$ are identically distributed (and $\log{\zeta_{n,k}(1)}$, $\ldots, \log{\zeta_{n,k}}(n)$ are identically distributed as well). Hence \eqref{a1} yields
\begin{align}
    \begin{gathered}
        \var\big(\widehat{D}_{n,m}(k,l)\big)
        =\frac{1}{n^2} \sum_{i,j=1}^n \cov\Big( \log{\phi_{m,l}(i)} - \log{\zeta_{n,k}(i)},  \log{\phi_{m,l}(j)} - \log{\zeta_{n,k}(j)}\Big)  \\
        =\frac{1}{n} \var \left( \log{\phi_{m,l}(1)} \right) + \frac{2}{n^2} \sum_{1 \leq i < j \leq n} \cov \left( \log{\phi_{m,l}(i)}, \log{\phi_{m,l}(j)} \right)  \\+ \frac{1}{n} \var \left( \log{\zeta_{n,k}(1)} \right) + \frac{2}{n^2} \sum_{1 \leq i < j \leq n} \cov \left( \log{\zeta_{n,k}(i)}, \log{\zeta_{n,k}(j)} \right)  \\
        -\frac{2}{n^2} \sum_{i,j=1}^{n} \cov \left( \log{\phi_{m,l}(i)}, \log{\zeta_{n,k}(j)} \right).
    \end{gathered}
\end{align}

We do not strictly adhere to notation used in
Theorem~\ref{th1} proof. Namely, the choice of the sets $A \subset \mathbb{R}^d$, $\widetilde{A} \subset \mathbb{R}^d$, positive $U_j, C_j(x), \widetilde{C}_j(x)$ and
integers $m_j, n_j$, where
$j\in \mathbb{Z}_+$ and $x\in \mathbb{R}^d$, could be different.
The proof of Theorem \ref{th_main2} is also divided into several steps. \textit{Steps 1-3} are devoted to the demonstration of relation $\frac{1}{n} \var{(\log\phi_{m,l}(1))} \to 0$ as $n, m \to \infty$, while \textit{Step 4} contains the proof of relation $\frac{2}{n^2} \sum_{1 \leq i < j \leq n} \cov(\log\phi_{m,l}(i), \log\phi_{m,l}(j))\to 0$ as $n,m\to \infty$.
In \textit{Step 5} we establish that
$$
\frac{2}{n^2} \sum_{1 \leq i < j \leq n} \cov(\log\zeta_{n,k}(i), \log\zeta_{n,k}(j))\to 0,\;\;n \to \infty,
$$
This step is rather involved.
In \textit{Step 6} we come to the desired statement $\var \left(\widehat{D}_{n,m}(k,l)\right) \to 0$, $n,m\to \infty$.

\textit{Step 1.} We  study $\e \log^2\left( \phi_{m,l}(1) \right)$, as $m\to \infty$. Consider
\begin{equation}\label{A}
A :=\Lambda(q) \cap S(q) \cap D_q(R_2)\cap S(p) \cap A_{p,2}(G_{N_1}),
\end{equation}
where the first four sets appeared in Theorem \ref{th1} proof, and $A_{p,2}(G_{N})$, for   $N \in \mathbb{N}$ and a probability density $p$ on $\mathbb{R}^d$, is defined quite similar to $A_{p}(G_{N})$. Namely, for $x\in \mathbb{R}^d$ and $N \in \mathbb{N}$, introduce
\begin{equation}\label{R2}
R_{N,2}(x) := \int_{\norm{x-y} \geq e_{[N]}} G_{N}(\log^2{\norm{x-y}}) q(y) \, dy
\end{equation}
and set $A_{p,2}(G_N):=\{x\in S(p): R_{N,2}(x) < \infty\}$.
Then $\p_X(S(p) \setminus A_{p,2}(G_{N_1})) = 0$
since $K_{p,q}(2, N_1) < \infty$. It is easily seen that $\p_{X}(A) = 1$. The reasoning is the same as in the proof of Theorem \ref{th1}.

Recall that, for each $x \in A$, one has $\log \xi_{m,l,x}\stackrel{law}\rightarrow \log \xi_{l,x},\,m\to \infty$, where $\xi_{m,l,x} := m \norm{x - Y_{(l)} (x,\mathbb{Y}_m)}^d$ and $\xi_{l,x}$ has $\Gamma(V_d\,q(x),l)$ distribution. Convergence in law of random variables is preserved under continuous mapping. Hence, for any
$x\in  A$, we come to the relation
\begin{equation}\label{b2**2}
\log^2 \xi_{m,l,x}\stackrel{law}\rightarrow \log^2 \xi_{l,x},\;\;m\to \infty.
\end{equation}
In view of \eqref{eq1a}, for each $x \in A$,
\begin{align}\label{mega_obvious}
    \begin{gathered}
        \e \log^2{\xi_{m,l,x}} = \int_{(0,\infty)} \log^2{u} \, dF_{m,l,x}(u) = \int_{(0,\infty)} \log^2{u} \, d
        \p(\phi_{m,l}(1) \leq u|X_1 = x)  \\=
        \e(\log^2{\phi_{m,l}(1)} | X_1 = x).
    \end{gathered}
\end{align}
Note that if $\eta\sim \Gamma(\alpha, \lambda)$, where $\alpha > 0$ and $\lambda >0$, then
\begin{align*}
    \begin{gathered}
        \e \log^2 \eta = \int_{(0,\infty)} \log^2{u} \, \frac{\alpha^\lambda u^{\lambda-1} e^{-\alpha u}}{\Gamma(\lambda)} \,du =
        \frac{1}{\Gamma{(\lambda)}} \int_{(0, \infty)} \left( \log{\frac{v}{\alpha}} \right)^2 \, v^{\lambda-1} e^{-v} \, dv\\ = \frac{1}{\Gamma{(\lambda)}}  \left(\int_{(0, \infty)} v^{\lambda-1} \log^2{v} \, e^{-v} \, dv - 2 \log{\alpha} \int_{(0, \infty)} v^{\lambda-1} \log{v} \, e^{-v} \, dv + \log^2\alpha\int_{(0, \infty)} v^{\lambda-1} e^{-v} \, dv \right) \\
        = \frac{\Gamma''(\lambda) - 2 \log\alpha \, \Gamma'(\lambda) + \log^2 \alpha \,  \Gamma(\lambda)}{\Gamma(\lambda)} = \frac{\Gamma''(\lambda)}{\Gamma(\lambda)} - 2 \, \psi(\lambda) \, \log\alpha + \log^2\alpha.
    \end{gathered}
\end{align*}
Since $\xi_{l,x} \sim \Gamma(V_d q(x), l)$ for $x \in S(q)$, one has
\begin{align}\label{calculate_integral_erlang_squared}
    \begin{gathered}
        \e \log^2{\xi_{l,x}} = \frac{\Gamma''(l)}{\Gamma(l)} - 2 \, \psi(l) \, \log(V_d q(x)) + \log^2(V_d q(x)) \\= \log^2 q(x) + \log q(x) \, \big( 2 \log V_d - 2 \psi(l) \big) + \left(\log^2 V_d - 2 \psi(l) \log V_d + \frac{\Gamma''(l)}{\Gamma(l)}\right)  \\=\log^2 q(x) + h_1 \log q(x) + h_2,
    \end{gathered}
\end{align}
where $h_1 := h_1(l,d)$ and $h_2 := h_2(l,d)$ depends only on fixed $l$ and $d$.

We  prove now that, for $x\in A$, one has
\begin{equation}\label{eauxil}
\e(\log^2{\phi_{m,l}(1)} | X_1 \!=\! x) \!\to\! \log^2 q(x)\! +\! h_1 \log q(x) \!+ \!h_2,\;\;m \to \infty.
\end{equation}
By virtue  of \eqref{mega_obvious} and \eqref{calculate_integral_erlang_squared}
relation \eqref{eauxil} is equivalent to the following one  $\e \log^2{\xi_{m,l,x}} \to
\e \log^2{\xi_{l,x}}$, $m \to \infty$.
So, in view of \eqref{b2**2} to prove \eqref{eauxil} it is sufficient to show that, for each $x \in A$, a family $\left\{\log^2 \xi_{m,l,x}\right\}_{m \geq m_0(x)}$ is uniformly integrable for some $m_0(x) \in \mathbb{N}$.  As in the proof of Theorem~\ref{th1}, we can verify that, for all $x \in A$ and some nonnegative $C_0(x)$,
\begin{equation}\label{b3_squared}
\sup_{m\geq m_0(x)} \e G_{N_1}(\log^2 \xi_{m,l,x})\leq C_0(x)<\infty.
\end{equation}

{\it Step 2}. Now our goal is to prove \eqref{b3_squared}.
For each $N \in \mathbb{N}$, introduce $\rho(N):=\exp\{\sqrt{e_{[N-1]}}\}$ and
$$
h_N(t):=
\begin{cases}
0, & t \in  \left(\frac{1}{\rho(N)}, \rho(N)\right],\\
\frac{2 \log{t}}{t}\left(\log_{[N]}(\log^2{t}) +  \frac{1}{\prod_{j=1}^{N-1} \log_{[j]}(\log^2{t})}\right),
& t \in \left(0, \frac{1}{\rho(N)}\right] \cup \left(\rho(N), \infty\right).
\end{cases}
$$
As usual, a product over an empty set (if $N=1$) is equal to $1$.

To show \eqref{b3_squared} we employ the following result.

\begin{lem}\label{lemma_G2}
Let $F(u), u\in \mathbb{R}$, be a distribution function such that $F(0)=0$. Fix an arbitrary  $N \in \mathbb{N}$. Then

1) $\int_{\left(0,\frac{1}{\rho(N)}\right]} G_N(\log^2 u) dF(u) = \int_{\left(0,\frac{1}{\rho(N)}\right]}F(u)(-h_N(u)) du$,

2) $\int_{\left(\rho(N), \infty\right)}G_N(\log^2 u) dF(u) = \int_{\left(\rho(N), \infty\right)}(1-F(u)) h_N(u) du$.
\end{lem}

The proof of this lemma is omitted, being quite similar to one of Lemma \ref{lemma_G}. By Lemma \ref{lemma_G2} and since $G_{N_1}(\log^2 u) = 0$, for $u \in \left( \frac{1}{\rho(N_1)}, \rho(N_1) \right]$, one has
\begin{gather*}
    \e G_{N_1}(\log^2{\xi_{m,l,x}}) = \int_{\left(0,\frac{1}{\rho(N_1)}\right]} F_{m,l,x}(u)(-h_{N_1}(u))du +
    \int_{\left(\rho(N_1), \infty\right)}(1-F_{m,l,x}(u))h_{N_1}(u)du \\
    :=I_1(m,x)+I_2(m,x).
\end{gather*}
To simplify notation we do not indicate the dependence of $I_i(m,x)$ ($i=1,2$) on $N_1$, $l$ and $d$.

We divide further proof into several parts.

{\it Part (2a)}. At first we consider $I_1(m,x)$.
As in Theorem \ref{th1} proof, for fixed $R_1 > 0$ and $\varepsilon_1 > 0$ appearing in the conditions of Theorem \ref{th_main2}, an inequality
   $ F_{m,l,x}(u) \leq (M_{q}(x, R_1))^{\varepsilon_1} V_d^{\varepsilon_1} u^{\varepsilon_1} $
holds, for any $x\in A$, $u \in \left( 0, \frac{1}{\rho(N_1)} \right]$ and $m \geq m_1 := \max\left\{ \left\lceil \frac{1}{\rho{(N_1)} R_1^d} \right\rceil, l \right\}$.
Taking into account that $0\leq (-h_{N_1}(u)) \leq \frac{(-2 \log u) \left(\log_{[N_1]}(\log^2 u) +1\right)}{u}$ if $u \in \left( 0, \frac{1}{\rho(N_1)} \right]$, we get, for $m \geq m_1$,
\begin{align}\label{third_proof_l2: eq2}
\begin{gathered}
I_1(m,x) \leq
(M_q(x, R_1))^{\varepsilon_1} V_d^{\varepsilon_1} \int_{\left(0,\frac{1}{\rho(N_1)}\right]}\!\!
\frac{(-2 \log u) \left(\log_{[N_1]}(\log^2 u)\! +\!1\right)}{u^{1-\varepsilon_1}} du\\
= U_1(\varepsilon_1, N_1, d) (M_q(x, R_1))^{\varepsilon_1}.
\end{gathered}
\end{align}
Here $U_1(\varepsilon, N, d) := V_d^{\varepsilon} L_{N,2}(\varepsilon)$, $L_{N,2}(\varepsilon) := \int_{\left[ \sqrt{e_{[N-1]}}, \infty \right)}
2 t \left( \log_{[N]}(t^2) + 1 \right) e^{-\varepsilon t} \, dt < \infty$
for each $\varepsilon > 0$ and any  $N \in \mathbb{N}$.

{\it Part (2b)}. Consider $I_2(m,x)$. As in the proof of Theorem \ref{th1}, taking into account that,
for $u \in (\rho(N_1), \infty)$,  $h_{N_1}(u) \leq \frac{2 \log{u}}{u} \left( \log_{[N_1]}(\log^2 u) + 1 \right)$,  we write, for all $m \geq \max\{\rho^2(N_1), l\}$,
\begin{align*}
\begin{gathered}
I_2(m,x) \leq
    \int_{(\rho(N_1), {\sqrt{m}}]} (1 - F_{m,l,x}(u)) \frac{2 \log{u} \, \left( \log_{[N_1]}(\log^2 u) + 1 \right)}{u} \, du \,\,  \\ + \int_{({\sqrt{m}},m^2]} (1 - F_{m,l,x}(u)) \frac{2 \log{u} \, \left( \log_{[N_1]}(\log^2 u) + 1 \right)}{u} \, du  \\ +\int_{(m^2, \infty)} (1 - F_{m,l,x}(u)) h_{N_1}(u) \, du  := J_1(m,x) + J_2(m,x) + J_3(m,x),
\end{gathered}
\end{align*}
where we do not indicate the dependence of $J_j(m,x)$ ($j=1,2,3$) on $N_1$ and $l$.

For $R_2 > 0$ and $\varepsilon_2 > 0$ appearing in the conditions of  Theorem \ref{th_main2}, one can prove (see Theorem \ref{th1} proof), that inequality
\begin{equation}\label{Fest}
1 - F_{m,l,x}(u)\leq S_1 \left( S_2 V_d u \, m_q(x, R_2) \right)^{-\varepsilon_2}
\end{equation}
holds for any $x \in A$, $u \in \left(\rho(N_1), \sqrt{m}\right]$ and all $m \geq m_2 := \max\left\{ \left\lceil \frac{1}{R_2^{2d}} \right\rceil, \left\lceil \rho^2(N_1) \right\rceil, l \right\}$. Here $S_1 := S_1(l)$ and $S_2$ are the same as in the proof of  Theorem \ref{th1}. For all $x \in A$ and $m \geq m_2$,  we come to the relations
\begin{align}\label{third_proof_l2: eq3}
\begin{gathered}
    J_1(m,x)
     \leq \frac{S_1}{(S_2 \, V_d)^{\varepsilon_2}(m_q(x, R_2))^{\varepsilon_2}}
     \int_{(\rho(N_1), \infty)}
     \frac{2 \log{u} \, \left( \log_{[N_1]}(\log^2 u) + 1 \right)}{u^{1+\varepsilon_2}} \, du\\
= U_2(\varepsilon_2, N_1, d, l) (m_q(x, R_2))^{-\varepsilon_2},
\end{gathered}
\end{align}
where $U_2(\varepsilon, N, d, l) := 2 S_1(l) \, L_{N,2}(\varepsilon)(S_2 \, V_d)^{-\varepsilon_2}$.

{\it Part (2c)}. Now we turn to $J_2(m,x)$.
Take  $\delta > 0$.
Then, due to \eqref{Fest}, for all $x \in A$ and any
$m \geq m_2$,
\begin{align}\label{U3_proof_l2_inequality}
    \begin{gathered}
        J_2(m,x)
        \leq  2 \left(1-F_{m,l,x}(\sqrt{m})\right) \int_{\left(\sqrt{m}, \, m^2\right]} \log{u} \,
        \left(\log_{[N_1]}(\log^2 u) + 1\right)\, d \log{u}
        \\ \leq 4 S_1 (S_2 V_d)^{-\varepsilon_2} m^{-\frac{\varepsilon_2}{2}}\left(m_q(x, R_2)\right)^{-\varepsilon_2}  \left( \log_{[N_1]}(4 \log^2 m) + 1 \right) \log^2 m \\
        = U_3(m, \varepsilon_2, N_1, d, l) \left(m_q(x, R_2)\right)^{-\varepsilon_2},
    \end{gathered}
\end{align}
where $U_3(m, \varepsilon, N, d, l) :=  4 S_1 (S_2 V_d)^{-\varepsilon_2} m^{-\frac{\varepsilon_2}{2}} \left( \log^2 m\right) \left( \log_{[N_1]}(4 \log^2 m) + 1 \right)\to 0$, $m \to \infty$.
\vskip0.2cm
{\it Part (2d)}. Now we consider $J_3(m,x)$. Take $u = m w$. Then $J_3(m,x)$ has the form
\begin{align*}
    \begin{gathered}
    \int_{\left(m, \, \infty \right)} (1-F_{m,l,x}(mw)) \frac{2 \log{(m w)}}{w}\left(\log_{[N_1]}(\log^2 (m w)) +  \frac{1}{\prod_{j=1}^{N_1-1} \log_{[j]}(\log^2 (m w))}\right) \, dw.
          \end{gathered}
\end{align*}
Due to Lemma \ref{sv} there exists $T(N)> \rho(N)$ such that
\begin{equation}\label{svf}
\log_{[N]}(\log^2(w^2))=\log_{[N]}(4 \log^2 w) \leq 2 \log_{[N]}(\log^2 w),\;\;w \geq T(N).
\end{equation}
Pick some $\delta > 0$ and set $m_3:=\max\left\{ l, \left\lceil (l-1)
\left( 1 + \frac{1}{\delta}\right) \right\rceil, \left\lceil T(N_1) \right\rceil, \left\lceil \rho(N_1) \right\rceil \right\}$, where $T(N)$ was introduced in \eqref{svf}.
Consider $m \geq m_3$.
In view of Lemma \ref{l4} (for $N = N_1$ and $\nu = 2$), \eqref{another_way_to_estimate_with_mw}, \eqref{svf}, \eqref{K} and
since $w > m$,
\begin{align*}
    \begin{gathered}
    J_3(m,x)\leq  \int_{\left(m, \, \infty \right)} (1-F_{m,l,x}(m w)) \frac{2 \log{( w^2)}}{w}\left(\log_{[N_1]}(\log^2 (w^2)) +  \frac{1}{\prod_{j=1}^{N_1-1} \log_{[j]}(\log^2 w)}\right) \, dw  \\
      \leq 4 (1 + \delta) \int_{\left(m, \, \infty \right)} (1-F_{1,1,x}(w)) \frac{2 \log{w}}{w}\left(\log_{[N_1]}(\log^2 w) +  \frac{1}{\prod_{j=1}^{N_1-1} \log_{[j]}(\log^2 w)}\right) \, dw \\
      =4 (1 + \delta) \int_{\left(m, \, \infty \right)} (1-F_{1,1,x}(w)) h_{N_1}(w) \, dw \leq
      4 (1 + \delta) \int_{\left(\rho(N_1), \, \infty \right)} (1-F_{1,1,x}(w)) h_{N_1}(w) \, dw   \end{gathered}
       \end{align*}
       \begin{align}\label{J3_in_l2_proof}
       \begin{gathered}
            =4 (1 + \delta) \int_{\left(\rho(N_1), \infty \right)} G_{N_1}(\log^2{w}) \, d F_{1,1,x}(w) =
      4 (1 + \delta)\e [G_{N_1}(\log^2 \xi_{1,1,x}) \ind\{\xi_{1,1,x} > \rho(N_1) \} ]  \\
      =4 (1 + \delta)\e [G_{N_1}((\log \norm{Y-x}^d)^2) \ind\{\|Y-x\|^d > \rho(N_1) \} ]  \\
    =4 (1 + \delta)\int_{y\in \mathbb{R}^d, \|x-y\| > (\rho(N_1))^{1/d}} G_{N_1}((\log \|x-y\|^d)^2) q(y) \, dy  \\\leq
    4 (1 + \delta) \left( a(d,2) \int_{y\in \mathbb{R}^d, \|x-y\| > \left(\rho(N_1)\right)^{1/d}} G_{N_1}(\log^2 \|x-y\|) q(y) \, dy + b(N_1,d, 2) \right)  \\
   = 4 (1 + \delta) \left( a(d,2) \left(  R_{N_1, 2}(x) + G_{N_1}(e^2_{[N_1-1]}) \right) + b(N_1,d, 2) \right) \\
   = A(\delta, d) R_{N_1, 2}(x) + B(\delta, d, N_1),
    \end{gathered}
\end{align}
$R_{N,2}(x)$ is defined in \eqref{R2}, $A(\delta, d) := 4(1+\delta) a(d,2)$,
$B(\delta, d, N_1) := 4(1 + \delta) \big( a(d,2) G_{N_1}(e^2_{[N_1-1]}) + b(N_1, d, 2) \big)$.

{\it Part (2e)}. Thus, for each $x\in A$ and $m \geq \max\{ m_1, m_2, m_3 \}$, taking into account \eqref{third_proof_l2: eq2}, \eqref{third_proof_l2: eq3}, \eqref{U3_proof_l2_inequality} and \eqref{J3_in_l2_proof}, we can claim that
\begin{align}\label{final_inequality_in_l2_proof_first}
    \begin{gathered}
        \e G_{N_1}(\log^2{\xi_{m,l,x}}) \leq I_1(m,x) + J_1(m,x) + J_2(m,x) + J_3(m,x)  \\ \leq U_1(\varepsilon_1, N_1, d) (M_q(x, R_1))^{\varepsilon_1} + U_2(\varepsilon_2, N_1, d, l) (m_q(x, R_2))^{-\varepsilon_2} \\  + U_3(m, \varepsilon_2, N_1, d, l) \left(m_q(x, R_2)\right)^{-\varepsilon_2} + A(\delta, d) R_{N_1, 2}(x) + B(\delta, d, N_1).
    \end{gathered}
\end{align}
Moreover, for any $\kappa > 0$, one can choose $m_4 := m_4(\kappa, \varepsilon_2, N_1, d, l) \in \mathbb{N}$ such that, for $m \geq m_4$, it holds $U_3(m, \varepsilon_2, N_1, d, l) \leq \kappa$. Then by \eqref{final_inequality_in_l2_proof_first}, for each $x\in A$ and $m \geq m_0 := \max\{ m_1, m_2, m_3, m_4 \}$,
\begin{align}\label{final_inequality_in_l2_proof_first_kappa}
    \begin{gathered}
        \e G_{N_1}(\log^2{\xi_{m,l,x}}) \leq U_1(\varepsilon_1, N_1, d) (M_q(x, R_1))^{\varepsilon_1} \\ + \big(U_2(\varepsilon_2, N_1, d, l) + \kappa\big) (m_q(x, R_2))^{-\varepsilon_2}   +         A(\delta, d) R_{N_1, 2}(x) + B(\delta, d, N_1) := C_0(x) < \infty.
    \end{gathered}
\end{align}
Hence we have proved the uniform integrability of the family $\left\{\log^2{\xi_{m,l,x}}\right\}_{m \geq m_0}$ for each $x \in A$. Therefore, for any $x\in A$ (thus for $\p_{X}$-almost every $x\in S(p)$),
relation \eqref{eauxil} holds.

{\it Step 3}. Now we can return to $\e \log^2 \phi_{m,l}(1)$. Set $\Delta_{m,l}(x):=\e(\log^2\phi_{m,l}(1)|X_1=x) = \e \log^2{\xi_{m,l,x}}$.
Consider $x\in A$ and take any $m\geq m_0$.  Function $G_{N_1}$ is nondecreasing and convex according to Lemma \ref{l6}. Due to the Jensen inequality
\begin{equation}\label{fourth: eq1_proof_l2}
\begin{gathered}
    G_{N_1}(\Delta_{m,l}(x)) = G_{N_1}(\e \log^2{\xi_{m,l,x}})
     \leq  \e G_{N_1}(\log^2{\xi_{m,l,x}}).
\end{gathered}
\end{equation}
Relation \eqref{fourth: eq1_proof_l2} guarantees that, for each $x\in A$ and all $m \geq m_0$,
\begin{gather*}
    \int_{\mathbb{R}^d} G_{N_1}(\Delta_{m,l}(x)) p(x) \, dx \leq
    U_1(\varepsilon_1, N_1, d) Q_{p,q}(\varepsilon_1, R_1) +  \big(U_2(\varepsilon_2, N_1, d, l) + \kappa\big) T_{p,q}(\varepsilon_2, R_2)  \\+ A(\delta, d) K_{p,q}(2, N_1) + B(\delta, d, N_1)<\infty.
\end{gather*}
We have established  uniform integrability of  the family $\{\Delta_{m,l}(\cdot)\}_{m \geq m_0}$
(w.r.t. measure $\p_{X}$). Therefore, we conclude that
\begin{gather*}
\e \log^2 \phi_{m,l}(1) \to \int_{\mathbb{R}^d} p(x) \log^2 q(x) \, dx + h_1 \int_{\mathbb{R}^d} p(x) \log q(x) \, dx + h_2, \;\; m \to \infty.
\end{gather*}
It is easily seen that finiteness of integrals $Q_{p,q}(\varepsilon_1, R_1)$, $T_{p,q}(\varepsilon_2, R_2)$ implies that
$$
\int_{\mathbb{R}^d}p(x) \log^2 q(x) dx <\infty,\;\;\int_{\mathbb{R}^d} p(x) |\log q(x)| dx <\infty.
$$
This is verified as in Remark \ref{rem1add} by taking into account that $\log^2 z \leq \frac{4}{\varepsilon^2} z^{\varepsilon}$ for all $z \geq 1$ and $\varepsilon > 0$.
Thus, $\e \log^2 \phi_{m,l}(1) \to \tau_2 < \infty$. Hence $\var\left( \log\phi_{m,l}(1) \right) = \e \log^2 \phi_{m,l}(1) - \left(\e \log \phi_{m,l}(1)\right)^2 \to \tau_2 - \tau^2_1  < \infty$, $m \to \infty$, where $\tau_1 := \psi(l) - \log{V_d} - \int_{\mathbb{R}^d} p(x) \log{q(x)} \, dx$ according to \eqref{mr}. Consequently, $\frac{1}{n} \var\left( \log\phi_{m,l}(1) \right) \to 0$ as $n, m \to \infty$.

\textit{Step 4.} Now we consider $\cov(\log\phi_{m,l}(i), \log\phi_{m,l}(j))$ for $i \neq j$, where $i,j\in \{1,\ldots,n\}$.
For $x,y\in \mathbb{R}^d$, introduce  conditional distribution function
\begin{equation}\label{mutual_cdf}
    \Phi^{i,j}_{m,l,x,y}(u,w) := \p(\phi_{m,l}(i) \leq u, \phi_{m,l}(j) \leq w | X_i = x, X_j = y), \;\;u,w\geq 0.
\end{equation}
For $x,y\in \mathbb{R}^d$, $u,w \geq 0$, $i\neq j$,
\begin{align}\label{t2: third: eq1_th2}
    \begin{gathered}
        \Phi^{i,j}_{m,l,x,y}(u,w) = 1 - \p(\phi_{m,l}(i) > u | X_i = x, X_j = y) \\ - \p(\phi_{m,l}(j) > w | X_i = x, X_j = y) + \p(\phi_{m,l}(i) > u, \phi_{m,l}(j) > w | X_i = x, X_j = y) \\ =
        1 - \p\left(\norm{x - Y_{(l)}(x, \mathbb{Y}_m)} > r_m(u)\right) -\p\left(\norm{y - Y_{(l)}(y, \mathbb{Y}_m)} > r_m(w)\right) \\ + \p\left(\norm{x - Y_{(l)}(x, \mathbb{Y}_m)} > r_m(u), \norm{y - Y_{(l)}(y, \mathbb{Y}_m)} > r_m(w)\right).
    \end{gathered}
\end{align}
Here $r_m(a) = \left(\frac{a}{m}\right)^{\frac{1}{d}}$ for all $a \geq 0$, as previously. One can write $\Phi_{m,l,x,y}(u,w)$ instead of ${\Phi}^{i,j}_{m,l,x,y}(u,w)$, because the right-hand side of \eqref{t2: third: eq1_th2} does not depend on $i$ and $j$.

Set $A_1 := \big\{(x,y): x \in A, \, y \in A, \, x \neq y\big\}$ and $A_2 := \big\{(x,y): x \in A, \, y \in A, \, x = y\big\}$, where $A$ is introduced in \eqref{A}. Evidently, $\left(\p_{X}\otimes \p_{X}\right)(A_1) = 1$ and $\left(\p_{X}\otimes \p_{X}\right)(A_2) = 0$.

Consider $(x,y) \in A_1$. Obviously, for  any $a>0$,
$r_m(a)\to 0$,  as $m\to \infty$. For $(x,y)\in A_1$
we take $m_5=m_5(u,w,\norm{x-y}) := \left\lceil\left(\frac{2}{\norm{x-y}}\right)^d  \max\left\{ u, w \right\}\right\rceil$. Then $r_m(u) < \frac{\norm{x-y}}{2}$ and $r_m(w)<\frac{\norm{x-y}}{2}$ for all $m \geq m_5$.
Thus $B(x,r_m(u)) \cap B(y, r_m(w)) = \varnothing$ if $m\geq m_5$. Consequently, for $m \geq m_6(u,w,\norm{x-y}) := \max\big\{ m_5, 2(l-1)\big\}$,
\begin{align}\label{mutual_proba_th2}
\begin{gathered}
    \p\left(\norm{x - Y_{(l)}(x, \mathbb{Y}_m)} > r_m(u), \norm{y - Y_{(l)}(y, \mathbb{Y}_m)} > r_m(w)\right)  \\=
    \sum_{s_1=0}^{l-1} \sum_{s_2=0}^{l-1} \frac{m!}{s_1! s_2! (m-s_1-s_2)!} \left( W_{m,x}(u) \right)^{s_1} \left( W_{m,y}(w) \right)^{s_2} \left( 1 - W_{m,x}(u) - W_{m,y}(w) \right)^{m-s_1-s_2}.
\end{gathered}
\end{align}
In view of \eqref{eq1a}, \eqref{t2: third: eq1_th2} and \eqref{mutual_proba_th2}, one has
for $\Phi_{m,l,x,y}(u,w)$ the following representation
\begin{align}\label{f_mutual_in_proof_of_theorem2}
    \begin{gathered}
       1 - \sum_{s_1=0}^{l-1}\binom{m}{s_1} \left( W_{m,x}(u) \right)^{s_1} \left( 1 - W_{m,x}(u) \right)^{m-s_1} - \sum_{s_2=0}^{l-1}\binom{m}{s_2} \left( W_{m,y}(w) \right)^{s_2} \left( 1 - W_{m,y}(w) \right)^{m-s_2} \\ +
        \sum_{s_1=0}^{l-1} \sum_{s_2=0}^{l-1} \frac{m!}{s_1! s_2! (m-s_1-s_2)!} \left( W_{m,x}(u) \right)^{s_1} \left( W_{m,y}(w) \right)^{s_2} \left( 1 - W_{m,x}(u) - W_{m,y}(w) \right)^{m-s_1-s_2}.
    \end{gathered}
\end{align}
For any fixed $(x,y)\in A_1$ and $u, w > 0$,
\begin{align}\label{mutual_proba_th2_limit_step3}
\begin{gathered}
\frac{m!}{s_1! s_2! (m-s_1-s_2)!} \left( W_{m,x}(u) \right)^{s_1} \left( W_{m,y}(w) \right)^{s_2}
\to  \frac{(V_d\, u\, q(x))^{s_1}}{s_1!} \frac{(V_d\, w\, q(y))^{s_2}}{s_2!}, \; m \to \infty, \\
\left( 1-W_{m,x}(u)-W_{m,y}(w) \right)^{m-s_1-s_2}
\to  e^{-V_d \big( u q(x) + w q(y) \big)}, \;\; m \to \infty.
\end{gathered}
\end{align}
Then, according to \eqref{f_mutual_in_proof_of_theorem2}, \eqref{conv} and \eqref{mutual_proba_th2_limit_step3}, for all fixed $u, w >0$, $(x,y) \in A_1$, one has
\begin{align*}
    \begin{gathered}
        \Phi_{m,l,x,y}(u,w)
        \to  1 - \sum_{s_1=0}^{l-1} \frac{(V_d u q(x))^{s_1}}{s_1!} e^{-V_d u q(x)} - \sum_{s_2=0}^{l-1} \frac{(V_d w q(y))^{s_2}}{s_2!} e^{-V_d w q(y)} \\+ \sum_{s_1=0}^{l-1} \sum_{s_2=0}^{l-1} \frac{(V_d\, u\, q(x))^{s_1}}{s_1!} \frac{(V_d\, w\, q(y))^{s_2}}{s_2!} e^{-V_d \big( u q(x) + w q(y) \big)}  \\=
        \Big( 1 - \sum_{s_1=0}^{l-1} \frac{(V_d u q(x))^{s_1}}{s_1!} e^{-V_d u q(x)} \Big)  \Big( 1 - \sum_{s_2=0}^{l-1} \frac{(V_d w q(y))^{s_2}}{s_2!} e^{-V_d w q(y)} \Big)  \\= F_{l,x}(u) F_{l,y}(w) := \Phi_{l,x,y}(u,w),\;\;m \to \infty.
    \end{gathered}
\end{align*}
Thus $\Phi_{l, x,y}(\cdot,\cdot)$ is a distribution function of a vector
$\eta_{l,x,y}:=(\xi_{l,x},\xi_{l,y})$,
where $\xi_{l,x}\sim \Gamma(V_d q(x), l)$, $\xi_{l,y} \sim \Gamma(V_d q(y), l)$
and the components of $\eta_{l,x,y}$ are independent.
Observe also that $\Phi_{m,l,x,y}(\cdot,\cdot)$ is a distribution function of a random vector
$\eta_{m,l,x,y} : = (\xi_{m,l,x}, \xi_{m,l,y})$.

Consequently, we have shown that $\eta_{m,l,x,y}\stackrel{law}\rightarrow \eta_{l,x,y}$ as $m\to \infty$.
Therefore, for any $(x,y) \in A_1$,
$$
\log\xi_{m,l,x} \log\xi_{m,l,y} \stackrel{law}\rightarrow \log\xi_{l,x}\log\xi_{l,y},\;\;m\to \infty.
$$
Here we exclude a set of zero probability where random variables under consideration can be equal to zero.
Note that, for all $i,j \in \mathbb{N}$, $i\neq j$,
\begin{align}\label{ij_indep}
\begin{gathered}
\e  (\log\xi_{m,l,x} \log\xi_{m,l,y}) = \int_{(0,\infty)} \int_{(0, \infty)} \log{u} \log{w} \, d \Phi_{m,l,x,y}(u, w) \\ = \e  \big( \log\phi_{m,l}(i) \log\phi_{m,l}(j)| X_i = x, X_j = y \big).
\end{gathered}
\end{align}
Obviously, in view of \eqref{calculate_integral_erlang} and since $\xi_{l,x}$ and $\xi_{l,y}$ are independent, one has
$$\e  (\log{\xi_{l,x}} \log{\xi_{l,y}}) = \e  \log{\xi_{l,x}} \,\e \log{\xi_{l,y}} = (\psi{(l)} - \log{V_d} - \log{q(x)})  (\psi{(l)} - \log{V_d} - \log{q(y)}).
$$

Now we intend to verify that, for any $(x,y)\in A_1$,
\begin{align}\label{q1_in_kull_leib}
\begin{gathered}
    \e  \big( \log\phi_{m,l}(1) \log\phi_{m,l}(2) | X_1 = x, X_2 = y \big)   \\\to (\psi{(l)} - \log{V_d} - \log{q(x)})  (\psi{(l)} - \log{V_d} - \log{q(y)}),\;\;m \to \infty.
\end{gathered}
\end{align}
Equivalently, one can prove that, for each $(x,y)\in A_1$, $\e  (\log\xi_{m,l,x} \log\xi_{m,l,y}) \to \e  (\log{\xi_{l,x}} \log{\xi_{l,y}})$, $m \to \infty$.

{\it Part (4a)}. We  establish the uniform integrability of a family $\{\log\xi_{m,l,x} \log\xi_{m,l,y}\}_{m \geq m_0}$ for  $(x,y) \in A_1$. The function $G_{N_1}(\cdot)$ is nondecreasing and convex. Thus,
for any $(x,y)\in A_1$, following the proof of \textit{Step 2}, one can find $m_0$ (the same as in the proof of  \textit{Step 2}) such that, for all $m \geq m_0$,
\begin{align}\label{Step2: 1_kl_proof}
\begin{gathered}
    \e  G_{N_1}(|\log\xi_{m,l,x} \, \log\xi_{m,l,y}|)
    \leq \frac{1}{2} \left(\e G_{N_1}(\log^2\xi_{m,l,x}) + \e G_{N_1}(\log^2\xi_{m,l,y})\right)  \\
    \leq \frac{U_1}{2} \Big((M_q(x, R_1))^{\varepsilon_1} + (M_q(y, R_1))^{\varepsilon_1}\Big) + \frac{U_2 + \kappa}{2} \Big((m_q(x, R_2))^{-\varepsilon_2} + (m_q(y, R_2))^{-\varepsilon_2}\Big) \\
       + \frac{A}{2} \Big( R_{N_1, 2}(x) + R_{N_1, 2}(y) \Big) + B := \widetilde{C}_0(x,y).
\end{gathered}
\end{align}
Clearly, $U_1, U_2, \kappa, A, B$ do not depend on $x$ or $y$ by virtue of \eqref{final_inequality_in_l2_proof_first_kappa}.
Hence, for any $(x,y) \in A_1$, a family $\{\log\xi_{m,l,x} \log\xi_{m,l,y}\}_{m \geq m_0}$ is uniformly integrable. Therefore we come to \eqref{q1_in_kull_leib} for $(x,y)\in A_1$.

{\it Part (4b)}. Set $T_{m,l}(x,y) := \e  \big( \log\phi_{m,l}(1) \log\phi_{m,l}(2) | X_1 = x, X_2 = y \big)$ $= \e  (\log \xi_{m,l,x} \, \log \xi_{m,l,y})$, where $(x,y)\in A_1$.
Then \eqref{q1_in_kull_leib} means that $T_{m,l}(x,y) \to (\psi{(l)} - \log{V_d} - \log{q(x)}) (\psi{(l)} - \log{V_d} - \log{q(y)})$  for any $(x,y)\in A_1$, as $m \to \infty$. Note that
\begin{align}\label{Step2: 4}
\begin{gathered}
    G_{N_1}(|T_{m,l}(x,y)|) = G_{N_1}(|\e  \log \xi_{m,l,x} \, \log \xi_{m,l,y}|) \\ \leq G_{N_1}(\e |\log \xi_{m,l,x} \, \log \xi_{m,l,y}| ) \leq \e G_{N_1}(|\log \xi_{m,l,x} \, \log \xi_{m,l,y}|).
\end{gathered}
\end{align}
Due to \eqref{Step2: 1_kl_proof} and \eqref{Step2: 4} one can conclude that, for all $m \geq m_0$, as $\left(\p_{X}\otimes \p_{X}\right)(A_1) = 1$,
\begin{gather*}
\int_{\mathbb{R}^d} \int_{\mathbb{R}^d} G_{N_1}(|T_{m,l}(x,y)|) p(x) p(y) \, dx \, dy
= \;\;\;\;\;\;\int\!\!\!\!\!\!\!\!\!\!\!\!\!\!\!\!\!\int\limits_{(x,y) \in A_1} G_{N_1}(|T_{m,l}(x,y)|) p(x) p(y) \, dx \, dy \\
\leq U_1 \int_{\mathbb{R}^{d}} M_q^{\varepsilon_1}(x, R_1) p(x) \, dx + \Big(U_2 + \kappa\Big) \int_{\mathbb{R}^{d}} m_q^{-\varepsilon_2}(x, R_2) p(x) \, dx  +  A \int_{\mathbb{R}^d} R_{N_1, 2}(x) p(x) \, dx + B  \\ = U_1 Q_{p,q}(\varepsilon_1, R_1) + (U_2 + \kappa) T_{p,q}(\varepsilon_2, R_2) + A K_{p,q}(2, N_1) + B < \infty.
\end{gather*}
Hence, for $(x,y)\in A_1$, a family
$\big\{ T_{m,l}(x,y) \big\}_{m \geq m_0}$ is uniformly integrable w.r.t. $\p_{X}\otimes \p_{X}$.
Consequently,
\begin{align}\label{final_conv}
\begin{gathered}
    \int_{\mathbb{R}^d} \int_{\mathbb{R}^d}\! T_{m,l}(x,y) p(x) p(y) \, dx \, dy \\
    \to \int_{\mathbb{R}^d} \int_{\mathbb{R}^d}\! (\psi{(l)} - \log{V_d} - \log{q(x)})  (\psi{(l)} - \log{V_d} - \log{q(y)}) p(x) p(y) \, dx \, dy, \;\; m \to \infty.
\end{gathered}
\end{align}
Thus
\vspace{-0.2cm}
\begin{equation}\label{final_cov_phi_1}
    \e  \log\phi_{m,l}(1) \log\phi_{m,l}(2) \to \left(\psi{(l)} - \log{V_d} - \int_{\mathbb{R}^{d}} \log{q(x)} p(x) \, dx\right)^2, \;\; m \to \infty.
\end{equation}
On the other hand, taking also into account \eqref{mr}, we come to the relation
\vspace{-0.2cm}
\begin{equation}\label{final_cov_phi_2}
    \e \log\phi_{m,l}(1) \e\log\phi_{m,l}(2) \to \left(\psi{(l)} - \log{V_d} - \int_{\mathbb{R}^{d}} \log{q(x)} p(x) \, dx\right)^2.
\end{equation}
Therefore \eqref{final_cov_phi_1} and \eqref{final_cov_phi_2} imply that
$$\frac{2}{n^2} \sum_{1 \leq i < j \leq n} \cov \left( \log{\phi_{m,l}(i)}, \log{\phi_{m,l}(j)} \right) =\frac{n-1}{n}\cov(\log\phi_{m,l}(1), \log\phi_{m,l}(2))\to 0, \; n, m \to \infty$$.

\vspace{-0.2cm}

\textit{Step 5.} Now we consider $\cov(\log\zeta_{n,k}(i), \log\zeta_{n,k}(j))$ for $i \neq j$, where $i,j\in \{1,\ldots,n\}$.
Similar to \textit{Step 4}, for $x,y\in \mathbb{R}^d$ and $u,w > 0$, introduce a conditional distribution function
\vspace{-0.1cm}
\begin{align}\label{mutual_cdf_for_zetas}
\begin{gathered}
    \widetilde{\Phi}^{i,j}_{n,k,x,y}(u,w) := \p(\zeta_{n,k}(i) \leq u, \zeta_{n,k}(j) \leq w | X_i = x, X_j = y) \\
    = \p \left(\norm{x-X_{(k)}(x, \{X_s\}_{s \neq i,j} \cup \{y\})} \leq r_{n-1}(u), \norm{y-X_{(k)}(y, \{X_s\}_{s \neq i,j} \cup \{x\})} \leq r_{n-1}(w)\right) \\ :=
    \p(\widetilde{\eta}_{n,k,x}^{\,y,i,j} \leq u, \widetilde{\eta}_{n,k,y}^{\,x,i,j} \leq w), \;\;u,w\geq 0,
\end{gathered}
\end{align}
where $\widetilde{\eta}_{n,k,x}^{\,y,i,j} := (n-1) \norm{x-X_{(k)}(x, \{X_s\}_{s \neq i,j} \cup \{y\})}^d$. We  write further $\widetilde{\Phi}_{n,k,x,y}(u,w)$, $\widetilde{\eta}_{n,k,x}^{\,y}$ and $\widetilde{\eta}_{n,k,y}^{\,x}$ instead of $\widetilde{\Phi}^{i,j}_{n,k,x,y}(u,w)$, $\widetilde{\eta}_{n,k,x}^{\,y,i,j}$, $\widetilde{\eta}_{n,k,y}^{\,x,i,j}$, respectively (since
$X_1,X_2,\ldots$ are i.i.d. random vectors). Moreover, $\widetilde{\Phi}_{n,k,x,y}(u,w)$ is the distribution function of a random vector $\widetilde{\eta}_{n,k,x,y} := (\widetilde{\eta}_{n,k,x}^{\,y}, \widetilde{\eta}_{n,k,y}^{\,x})$ and the regular conditional distribution function of a random vector
$(\zeta_{n,k}(i), \zeta_{n,k}(j))$ given $(X_i, X_j) = (x,y)$. One has
\vspace{-0.2cm}
\begin{align*}
    \begin{gathered}
		\widetilde{\Phi}_{n,k,x,y}(u,w) = 1 - \p \left(\norm{x-X_{(k)}(x, \{X_s\}_{s \neq i,j} \cup \{y\})} > r_{n-1}(u)\right) \\ - \p \left(\norm{y-X_{(k)}(y, \{X_s\}_{s \neq i,j} \cup \{x\})} > r_{n-1}(w)\right) \\ + \p \left(\norm{x-X_{(k)}(x, \{X_s\}_{s \neq i,j} \cup \{y\})} > r_{n-1}(u), \norm{y-X_{(k)}(y, \{X_s\}_{s \neq i,j} \cup \{x\})} > r_{n-1}(w)\right).
    \end{gathered}
\end{align*}

Introduce
\vspace{-0.3cm}
\begin{equation*}
\widetilde{A} :=\Lambda(p) \cap S(p) \cap D_p(R_4)\cap \widetilde{A}_{p,2}(G_{N_2}),
\end{equation*}
where the first three sets appeared in Theorem \ref{th1} proof (\textit{Step 5}), and $\widetilde{A}_{p,2}(G_{N})$, for   $N \in \mathbb{N}$ and a probability density $p$ on $\mathbb{R}^d$, is defined in full similarity to $\widetilde{A}_{p}(G_{N})$. Namely, introduce
\begin{equation*}
\widetilde{R}_{N,2}(x) := \int_{\norm{x-y} \geq e_{[N]}} G_{N}(\log^2{\norm{x-y}}) p(y) \, dy
\end{equation*}
and set $\widetilde{A}_{p,2}(G_N):=\{x\in S(p): \widetilde{R}_{N,2}(x) < \infty\}$.
Then $\p_X(S(p) \setminus \widetilde{A}_{p,2}(G_{N_2})) = 0$ since $K_{p,p}(2, N_2) < \infty$. It is easily seen that $\p_{X}(\widetilde{A}) = 1$.

Consider $\widetilde{A}_1 := \big\{(x,y): x \in \widetilde{A}, \, y \in \widetilde{A}, \, x \neq y\big\}$ and $\widetilde{A}_2 := \big\{(x,y): x \in \widetilde{A}, \, y \in \widetilde{A}, \, x = y\big\}$. Evidently, $\left(\p_{X}\otimes \p_{X}\right)(\widetilde{A}_1) = 1$ and $\left(\p_{X}\otimes \p_{X}\right)(\widetilde{A}_2) = 0$.
For any $a>0$,
$r_m(a)\to 0$, as $m\to \infty$. Hence, for $(x,y)\in \widetilde{A}_1$,
one can find $\widetilde{n}_5=\widetilde{n}_5(u,w,\norm{x-y}) = 1+ \left\lceil\left(\frac{2}{\norm{x-y}}\right)^d  \max\left\{ u, w \right\}\right\rceil$ such that $r_{n-1}(u) < \frac{\norm{x-y}}{2}$, $r_{n-1}(w)<\frac{\norm{x-y}}{2}$ if $n \geq \widetilde{n}_5$.
Then $B(x,r_{n-1}(u)) \cap B(y, r_{n-1}(w)) = \varnothing$ if $n\geq \widetilde{n}_5(u,w,\norm{x-y})$. Thus, for $n \geq \widetilde{n}_6 := \max\big\{ \widetilde{n}_5, 2k\big\}$, one has
\begin{equation*}
	\widetilde{\Phi}_{n,k,x,y}(u,w) = 1 - \sum_{s_1=0}^{k-1}\binom{n-2}{s_1} \left( V_{n-1,x}(u) \right)^{s_1} \left( 1 - V_{n-1,x}(u) \right)^{n-2-s_1}
\end{equation*}
\vspace{-0.4cm}
\begin{align}\label{Phi_tilde_equality}
    \begin{gathered}
  - \sum_{s_2=0}^{k-1}\binom{n-2}{s_2} \left( V_{n-1,y}(w) \right)^{s_2} \left( 1 - V_{n-1,y}(w) \right)^{n-2-s_2}
    \end{gathered}
\end{align}
\vspace{-0.3cm}
\begin{equation*}
 + \sum_{s_1=0}^{k-1} \sum_{s_2=0}^{k-1} \frac{(n-2)!}{s_1! s_2! (n-2-s_1-s_2)!} \left( V_{n-1,x}(u) \right)^{s_1} \left( V_{n-1,y}(w) \right)^{s_2} \left( 1 - V_{n-1,x}(u) - V_{n-1,y}(w) \right)^{n-2-s_1-s_2}.	
\end{equation*}
\vspace{0cm}

Therefore, for each fixed $(x,y) \in \widetilde{A}_1$, $u, w >0$, we get, as $n\to \infty$,
\begin{align*}
    \begin{gathered}
        \widetilde{\Phi}_{n,k,x,y}(u,w)
        \to  1 - \sum_{s_1=0}^{k-1} \frac{(V_d u\, p(x))^{s_1}}{s_1!} e^{-V_d u\, p(x)} - \sum_{s_2=0}^{k-1} \frac{(V_d w\, p(y))^{s_2}}{s_2!} e^{-V_d w\, p(y)} \\+ \sum_{s_1=0}^{k-1} \sum_{s_2=0}^{k-1} \frac{(V_d\, u\, p(x))^{s_1}}{s_1!} \frac{(V_d\, w\, p(y))^{s_2}}{s_2!} e^{-V_d \big( u\, p(x) + w\, p(y) \big)} \\  =
        \Big( 1 - \sum_{s_1=0}^{k-1} \frac{(V_d u\, p(x))^{s_1}}{s_1!} e^{-V_d u\, p(x)} \Big)  \Big( 1 - \sum_{s_2=0}^{k-1} \frac{(V_d w\, p(y))^{s_2}}{s_2!} e^{-V_d w\, p(y)} \Big)  = \widetilde{F}_{k,x}(u) \widetilde{F}_{k,y}(w) \\ := \widetilde{\Phi}_{k,x,y}(u,w).
    \end{gathered}
\end{align*}
Here $\widetilde{\Phi}_{k, x,y}(\cdot,\cdot)$ is the distribution function of a vector
$\widetilde{\eta}_{k,x,y}:=(\widetilde{\xi}_{k,x},\widetilde{\xi}_{k,y})$,
where $\widetilde{\xi}_{k,x}\sim \Gamma(V_d\, p(x), k)$, $\widetilde{\xi}_{k,y} \sim \Gamma(V_d\, p(y), k)$
and the components of $\widetilde{\eta}_{k,x,y}$ are independent.

Consequently, we have shown that $\widetilde{\eta}_{n,k,x,y}\stackrel{law}\rightarrow \widetilde{\eta}_{k,x,y}$ as $n\to \infty$.
Therefore, for any $(x,y) \in \widetilde{A}_1$,
\begin{equation}\label{weak_convergence_for_production_of_log_zetas}
\log\widetilde{\eta}_{n,k,x}^{\,y} \log\widetilde{\eta}_{n,k,y}^{\,x} \stackrel{law}\rightarrow \log\widetilde{\xi}_{k,x}\log\widetilde{\xi}_{k,y},\;\;n\to \infty.
\end{equation}
Here we exclude a set of zero probability where random variables under consideration can be equal to zero.
In a similar way to \eqref{ij_indep}, for $i,j\in \{1,\ldots,n\}$, $i\neq j$, we write
\vspace{-0.3cm}
\begin{align}\label{ij_indep_for_zetas}
\begin{gathered}
\e  \log\widetilde{\eta}_{n,k,x}^{\,y} \log\widetilde{\eta}_{n,k,y}^{\,x} = \int_{(0,\infty)} \int_{(0, \infty)} \log{u} \log{w} \, d \widetilde{\Phi}_{n,k,x,y}(u, w) \\ = \e  \big( \log\zeta_{n,k}(i) \log\zeta_{n,k}(j)| X_i = x, X_j = y \big).
\end{gathered}
\end{align}
Since $\widetilde{\xi}_{k,x}$ and $\widetilde{\xi}_{k,y}$ are independent, formula \eqref{calculate_integral_erlang} yields
$$\e  (\log{\widetilde{\xi}_{k,x}} \log{\widetilde{\xi}_{k,y}}) = \e  \log{\widetilde{\xi}_{k,x}} \,\e \log{\widetilde{\xi}_{k,y}} = (\psi{(k)} - \log{V_d} - \log{p(x)})  (\psi{(k)} - \log{V_d} - \log{p(y)}).
$$

For any fixed $M>0$, consider $\widetilde{A}_{1,M} := \big\{ (x,y) \in \widetilde{A}_1: \norm{x-y} > M \big\}$. Now our aim is to verify that, for each $(x,y)\in \widetilde{A}_{1,M}$,
\begin{align}\label{q1_in_kull_leib_for_zetas}
\begin{gathered}
    \e  \big( \log\zeta_{n,k}(1) \log\zeta_{n,k}(2) | X_1 = x, X_2 = y \big)   \\\to (\psi{(k)} - \log{V_d} - \log{p(x)})  (\psi{(k)} - \log{V_d} - \log{p(y)}),\;\;n \to \infty.
\end{gathered}
\end{align}
Equivalently, we can prove, for each  $(x,y)\in \widetilde{A}_{1,M}$, that
\begin{equation}\label{convergence_of_expect_for_zetas}
	\e  \log\widetilde{\eta}_{n,k,x}^{\,y} \log\widetilde{\eta}_{n,k,y}^{\,x} \to \e  \log{\widetilde{\xi}_{k,x}} \log{\widetilde{\xi}_{k,y}}, \;\; n \to \infty.
\end{equation}

The idea that we consider only $(x,y) \in \widetilde{A}_{1,M}$ is principle for the further proof.

{\it Part (5a)}. We will establish the uniform integrability of a family $\{\log\widetilde{\eta}_{n,k,x}^{\,y} \log\widetilde{\eta}_{n,k,y}^{\,x}\}_{n \geq \widetilde{n}_0}$ for $(x,y) \in \widetilde{A}_{1,M}$ and some $\widetilde{n}_0 \in \mathbb{N}$ which does not depend on $x, y$, but can depend on $M$. Then, due to \eqref{weak_convergence_for_production_of_log_zetas}, the relation \eqref{convergence_of_expect_for_zetas} would be valid for such $(x,y)$ as well.

As we have seen, the function $G_{N_2}(\cdot)$ is nondecreasing and convex. Hence
\begin{align}\label{Step2: 1_kl_proof_for_zetas}
\begin{gathered}
    \e  G_{N_2}(|\log{\widetilde\eta}_{n,k,x}^{\,y} \, \log{\widetilde\eta}_{n,k,y}^{\,x}|)
    \leq \frac{1}{2} \left(\e G_{N_2}(\log^2\widetilde\eta_{n,k,x}^{\,y}) + \e G_{N_2}(\log^2\widetilde\eta_{n,k,y}^{\,x})\right).
\end{gathered}
\end{align}

Let us consider, for instance, $\e G_{N_2} (\log^2 \widetilde\eta_{n,k,x}^{\,y})$. As at \textit{Step 2} we can write
\begin{align*}
	\begin{gathered}
		\e G_{N_2} (\log^2 \widetilde\eta_{n,k,x}^{\,y}) = \int_{\left(0,\frac{1}{\rho(N_2)}\right]} \widetilde{F}_{n,k,x}^{\,y}(u)(-h_{N_2}(u))du +
    	\int_{\left(\rho(N_2), \infty\right)}(1-\widetilde{F}_{n,k,x}^{\,y}(u))h_{N_2}(u)du \\ := I_1(n,x,y) + I_2(n,x,y),
	\end{gathered}
\end{align*}
where
\vspace{-0.3cm}
\begin{align}\label{cdf_eta_nkxy}
	\begin{gathered}
		\widetilde{F}_{n,k,x}^{y}(u) := \p \left( \widetilde{\eta}_{n,k,x}^{\,y} \leq u \right) = 1 - \p \left(\norm{x-X_{(k)}(x, \{X_s\}_{s \neq i,j} \cup \{y\})} > r_{n-1}(u)\right) \\ = \ind{\left\{ \norm{x-y} > r_{n-1}(u) \right\}} \left(1 - \sum_{s=0}^{k-1}\binom{n-2}{s} \left( V_{n-1,x}(u) \right)^{s} \left( 1 - V_{n-1,x}(u) \right)^{n-2-s}\right) \\ + \ind{\left\{ \norm{x-y} \leq r_{n-1}(u) \right\}} \left(1 - \sum_{s=0}^{k-2}\binom{n-2}{s} \left( V_{n-1,x}(u) \right)^{s} \left( 1 - V_{n-1,x}(u) \right)^{n-2-s}\right).
	\end{gathered}
\end{align}
As usual a sum over empty set is equal to $0$ (for $k = 1$).

If $u \in \left(0, \frac{1}{\rho{(N_2)}}\right]$, where $\rho(N) := \exp\{\sqrt{e_{[N-1]}}\}$
and $n \geq \widetilde{n}_1 := \left\lceil\frac{1}{\rho{(N_2)} M^d}\right\rceil + 1$, then $r_{n-1}(u) \leq M$. Thus $r_{n-1}(u) < \norm{x-y}$. In view of \eqref{cdf_eta_nkxy}, $\widetilde{F}_{n,k,x}^{y}(u) = 1 - \sum_{s=0}^{k-1}\binom{n-2}{s} \big( V_{n-1,x}(u) \big)^{s} \\ ( 1 -V_{n-1,x}(u) )^{n-2-s}$. Similarly to \eqref{ineq_B}, one has
\vspace{-0.2cm}
\begin{equation}
	\widetilde{F}_{n,k,x}^{\,y}(u) \leq \left( \frac{n-2}{n-1}\right)^{\varepsilon_3} \left( M_p(x,R_3) V_d u \right)^{\varepsilon_3} \leq \left( M_p(x,R_3) \right)^{\varepsilon_3} V_d^{\varepsilon_3} u^{\varepsilon_3}
\end{equation}
for all $(x,y) \in \widetilde{A}_{1,M}$, $u \in \left(0, \frac{1}{\rho(N_2)}\right]$, $n \geq \max\{\widetilde{n}_1(M), \widetilde{n}_2(R_3)\}$, where $\widetilde{n}_2(R_3) := \max\big\{ \left\lceil \frac{1}{\rho(N_2) R_3^d} \right\rceil + 1, k+1 \big\}$. Consequently, $I_1(n,x,y) \leq U_1(\varepsilon_3, N_2, d) \left( M_p(x, R_3)\right)^{\varepsilon_3}$ for all $(x,y) \in \widetilde{A}_{1,M}$ and $n \geq \max\left\{ \widetilde{n}_1(M), \widetilde{n}_2(R_3) \right\}$.
Moreover, for all $u > 0$, in view of \eqref{cdf_eta_nkxy} it holds
\begin{equation}
	1 - \widetilde{F}_{n,k,x}^{\,y}(u) \leq \sum_{s=0}^{k-1}\binom{n-2}{s} \left( V_{n-1,x}(u) \right)^{s} \left( 1 - V_{n-1,x}(u) \right)^{n-2-s}.
\end{equation}
The same reasoning as was used in Theorem \ref{th1} proof (\textit{Step 3, Part (3b)}) leads to the inequalities
\begin{align}
	\begin{gathered}
		1 - \widetilde{F}_{n,k,x}^{\,y}(u) \leq S_1(k) \left( 1 - S_2 \, V_{n-1,x}(u) \right)^{n-2} \leq
    	S_1 \exp\left\{  - S_2 \, (n-2)  V_{n-1,x}(u) \right\} \\ \leq S_1 \exp\left\{ -\frac{n-2}{n-1} \, S_2 V_d u \, m_p(x, R_4)\right\} \leq S_1 \left( \frac{S_2}{2} V_d u \, m_p(x, R_4) \right)^{-\varepsilon_4}
	\end{gathered}	
\end{align}
for all $n \geq \max\left\{\widetilde{n}_3(R_4), 3\right\}$. Then similarly to \eqref{final_inequality_in_l2_proof_first}, the relation
\vspace{-0.2cm}
\begin{align*}
\begin{gathered}
        \e G_{N_2}(\log^2{\widetilde\eta_{n,k,x}^{\,y}}) \leq U_1 (M_p(x, R_3))^{\varepsilon_3} + \big(\widetilde{U}_2+ \kappa\big) (m_p(x, R_4))^{-\varepsilon_4} + A \, \widetilde{R}_{N_1, 2}(x) + B := \widetilde{C}_0(x) < \infty
\end{gathered}
\end{align*}
is valid for all $(x,y) \in \widetilde{A}_{1,M}$ and $n \geq \widetilde{n}_0(M) := \max\left\{ \widetilde{n}_1, \widetilde{n}_2, \widetilde{n}_3, \widetilde{n}_4(\kappa), 3 \right\}$. Here $U_1, \widetilde{U}_2, \kappa, A, B$ do not depend on $x$ or $y$.
Thus, in view of \eqref{Step2: 1_kl_proof_for_zetas}, one has
\begin{align}\label{e_g_n2_for_zetas}
\begin{gathered}
    \e  G_{N_2}(|\log\widetilde\eta_{n,k,x}^{\,y} \, \log\widetilde\eta_{n,k,y}^{\,x}|)
    \leq \frac{U_1}{2} \Big((M_p(x, R_3))^{\varepsilon_3} + (M_p(y, R_3))^{\varepsilon_3}\Big) \\ + \frac{U_2 + \kappa}{2} \Big((m_p(x, R_4))^{-\varepsilon_4} + (m_p(y, R_4))^{-\varepsilon_4}\Big)
       + \frac{A}{2} \Big( \widetilde{R}_{N_2, 2}(x) + \widetilde{R}_{N_2, 2}(y) \Big) + B := \widetilde{C}_0(x,y).
\end{gathered}
\end{align}
Hence, for any $(x,y) \in \widetilde{A}_{1,M}$, a family $\{\log\widetilde\eta_{n,k,x}^{\,y} \log\widetilde\eta_{n,k,y}^{\,x}\}_{n \geq \widetilde{n}_0}$ is uniformly integrable. Thus we come to \eqref{q1_in_kull_leib_for_zetas} for $(x,y)\in \widetilde{A}_{1,M}$.

{\it Part (5b)}. Set $\widetilde{T}_{n,k}(x,y) := \e  \big( \log\zeta_{n,k}(1) \log\zeta_{n,k}(2) | X_1 = x, X_2 = y \big)$ $= \e  \log \widetilde\eta_{n,k,x}^{\,y} \, \log \widetilde\eta_{n,k,y}^{\,x}$ for all $(x,y)\in \widetilde{A}_1$.
Relation \eqref{q1_in_kull_leib_for_zetas} validity is equivalent to the following one: for any $(x,y)\in \widetilde{A}_{1,M}$, $\widetilde{T}_{n,k}(x,y) \to (\psi{(k)} - \log{V_d} - \log{p(x)}) (\psi{(k)} - \log{V_d} - \log{p(y)})$, as $n \to \infty$.
Now take any $(x,y) \in \widetilde{A}_1$. Then, for any fixed $M > 0$ and $(x,y) \in \widetilde{A}_1$, we have proved that
\begin{align}\label{Step2: asconv_for_zetas}
\begin{gathered}
\widetilde{T}_{n,k}(x,y) \ind\{\norm{x-y} > M\}   \\ \to (\psi{(k)} - \log{V_d} - \log{p(x)}) (\psi{(k)} - \log{V_d} - \log{p(y)}) \ind\{\norm{x-y}> M\}, \; n \to \infty.
\end{gathered}
\end{align}
Note that
\begin{align}\label{Step2: 4_for_zetas}
\begin{gathered}
    G_{N_2}(|\widetilde{T}_{n,k}(x,y)| \ind\{\norm{x-y} > M\}) \leq G_{N_2}(|\widetilde{T}_{n,k}(x,y)|) = G_{N_2}(|\e  \log \widetilde\eta_{n,k,x}^{\,y} \, \log \widetilde\eta_{n,k,y}^{\,x}|) \\ \leq G_{N_2}(\e |\log \widetilde\eta_{n,k,x}^{\,y} \, \log \widetilde\eta_{n,k,y}^{\,x}| ) \leq \e G_{N_2}(|\log \widetilde\eta_{n,k,x}^{\,y} \, \log \widetilde\eta_{n,k,y}^{\,x}|).
\end{gathered}
\end{align}
Due to \eqref{e_g_n2_for_zetas} and \eqref{Step2: 4_for_zetas} one can conclude that, for all $n \geq \widetilde{n}_0$,
\begin{gather*}
\int_{\mathbb{R}^d} \int_{\mathbb{R}^d} G_{N_2}(|\widetilde{T}_{n,k}(x,y)| \ind\{\norm{x-y} > M\}) p(x) p(y) \, dx \, dy \\ \leq U_1 \int_{\mathbb{R}^{d}} M_p^{\varepsilon_3}(x, R_3) p(x) \, dx + \Big(\widetilde{U}_2 + \kappa\Big) \int_{\mathbb{R}^{d}} m_p^{-\varepsilon_4}(x, R_4) p(x) \, dx  +  A \int_{\mathbb{R}^d} \widetilde{R}_{N_1, 2}(x) p(x) \, dx + B  \\ = U_1 Q_{p,p}(\varepsilon_3, R_3) + (\widetilde{U}_2 + \kappa) T_{p,p}(\varepsilon_4, R_4) + A K_{p,p}(2, N_2) + B < \infty.
\end{gather*}
Hence, for $(x,y)\in \widetilde{A}_1$, a family
$\big\{ \widetilde{T}_{n,k}(x,y) \ind\{ \norm{x-y} > M\}\big\}_{n \geq \widetilde{n}_0}$ is uniformly integrable w.r.t. $\p_{X}\otimes \p_{X}$.
Consequently, in view of \eqref{q1_in_kull_leib_for_zetas}, for each $M>0$,
\begin{align}\label{final_conv_for_zetas}
\begin{gathered}
	\int\!\!\!\!\!\!\!\!\!\!\!\!\!\!\!\!\!\!\int\limits_{x,y\in \mathbb{R}^d,\norm{x-y} > M}\!\!\!\! \widetilde{T}_{n,k}(x,y) p(x) p(y) \, dx \, dy \\
		\to \int\!\!\!\!\!\!\!\!\!\!\!\!\!\!\!\!\!\!\int\limits_{x,y\in \mathbb{R}^d,\norm{x-y} > M}\!\!\!\! (\psi{(k)} - \log{V_d} - \log{p(x)})  (\psi{(k)} - \log{V_d} - \log{p(y)}) p(x) p(y) \, dx \, dy, \;\; n \to \infty.
\end{gathered}
\end{align}
Now we consider the case $\norm{x-y} \leq M$.
One has $\bigcap_{s=1}^{\infty} \left\{ \norm{X_1-X_2} \leq \frac{1}{s} \right\} = \left\{ X_1 = X_2 \right\}$ and $\p \left( X_1 = X_2 \right) = 0$ as $X_1$ and $X_2$ are independent and have a density $p(x)$ w.r.t. the Lebesgue measure $\mu$. Then
\vspace{-0.1cm}
\begin{equation*}
\p\Big(\norm{X_1-X_2}\leq M\Big) \to 0, \;\; M \to 0,
\end{equation*}
Taking into account that, for an integrable function $h$,  $\int_C h d\p\to 0$ as $\p(C)\to 0$, 
we get
\begin{equation}
	\e (\log\zeta_{n,k}(1) \log\zeta_{n,k}(2) \ind\{\norm{X_1-X_2} \leq M\}) \to 0, M \to 0,
\end{equation}
since $\e \log\zeta_{n,k}(1) \log\zeta_{n,k}(2) \leq \frac{1}{2} \left( \e \log^2\zeta_{n,k}(1) + \e \log^2\zeta_{n,k}(2) \right) < \infty$ (the proof is similar to the establishing that $\e \log\phi_{m,l}(1) < \infty$).
Hence, for any $\gamma > 0$, one can find $M_1=M_1(\gamma) > 0$ such that, for all $M \in (0,M_1]$ and $n \geq \widetilde{n}_0$,
\begin{gather*}
    \Bigg|\;\;\;\;\;\;\;\int\!\!\!\!\!\!\!\!\!\!\!\!\!\!\!\!\!\!\int\limits_{x,y\in \mathbb{R}^d,\,\norm{x-y} \leq M} \widetilde{T}_{n,k}(x,y) p(x) p(y) \, dx \, dy \Bigg| = \big|\e \log\phi_{m,l}(1) \log\phi_{m,l}(2) \ind\{\norm{X_1-X_2} \leq M\} \big| < \frac{\gamma}{3}.
\end{gather*}
Set $v(t):= \psi{(k)} - \log{V_d} - \log{p(t)}$, $t\in \mathbb{R}^d$.
Also
there exists $M_2=M_2(\gamma)>0$ such that, for all $M\in (0,M_2]$,
\begin{gather*}
\Bigg|\;\;\;\;\;\;\;\int\!\!\!\!\!\!\!\!\!\!\!\!\!\!\!\!\!\!\int\limits_{x,y\in \mathbb{R}^d,\,\norm{x-y} \leq M}v(x)v(y) p(x) p(y) \, dx  \, dy \Bigg| < \frac{\gamma}{3}.
\end{gather*}
Take $M= \min\{M_1,M_2\}$. Due to \eqref{final_conv_for_zetas} one can find $\widetilde{n}_7(M, \gamma)$ such that for all $n \geq \max\{ \widetilde{n}_0, \widetilde{n}_7(M, \gamma)\}$ the following inequality holds
\begin{gather*}
    \Bigg|\;\;\;\;\;\;\;\int\!\!\!\!\!\!\!\!\!\!\!\!\!\!\!\!\!\!\int\limits_{x,y\in \mathbb{R}^d,\,\norm{x-y} > M} \widetilde{T}_{n,k}(x,y) p(x) p(y) \, dx \, dy   -\;\;\;\;\;\;\;\;\int\!\!\!\!\!\!\!\!\!\!\!\!\!\!\!\!\!\!\int\limits_{x,y\in \mathbb{R}^d,\,\norm{x-y} > M}
v(x)v(y) p(x) p(y) \, dx \, dy \Bigg| < \frac{\gamma}{3}.
\end{gather*}
So, for any $\gamma > 0$, there is $M(\gamma) >0$ such that, for all $n \geq \max\{ \widetilde{n}_0, \widetilde{n}_7(M, \gamma)\}$, one has
\begin{gather}\label{final_joint_for_zetas}
  \bigg|\int_{\mathbb{R}^{d}} \int_{\mathbb{R}^{d}} \widetilde{T}_{n,k}(x,y) p(x) p(y) \, dx \, dy  -
\int_{\mathbb{R}^{d}} \int_{\mathbb{R}^{d}} v(x)v(y) p(x) p(y) \, dx \, dy \bigg|
< \gamma.
\end{gather}
By virtue of the formula
\begin{gather*}
\int_{\mathbb{R}^{d}} \int_{\mathbb{R}^{d}} v(x)v(y) p(x) p(y) \, dx \, dy =  \left(\psi{(l)} - \log{V_d} - \int_{\mathbb{R}^{d}} (\log{p(x)}) p(x) \, dx\right)^2,
\end{gather*}
and taking into account \eqref{final_joint_for_zetas} we come to the relation
\begin{gather*}
\e \log\zeta_{n,k}(1) \log\zeta_{n,k}(2) \to \left(\psi{(k)} - \log{V_d} - \int_{\mathbb{R}^{d}} (\log{p(x)}) p(x) \, dx\right)^2.
\end{gather*}
Moreover, in view of \eqref{mra} (see \textit{Step 5} of Theorem \ref{th1} proof), we have
\begin{gather*}
\e \log\zeta_{n,k}(1) \e \log\zeta_{n,k}(2) \to \left(\psi{(k)} - \log{V_d} - \int_{\mathbb{R}^{d}} (\log{p(x)}) p(x) \, dx\right)^2.
\end{gather*}
Therefore 
$$\frac{2}{n^2} \sum_{1 \leq i < j \leq n} \cov \left( \log{\zeta_{n,k}(i)}, \log{\zeta_{n,k}(j)} \right) = \frac{2(n-1)}{n}\cov(\log\zeta_{n,k}(1), \log\zeta_{n,k}(2))\to 0, \; n \to \infty.$$

\vspace{-0.4cm}

\textit{Step 6.} Reasoning as at \textit{Steps 1-3} shows that $\frac{1}{n} \var \left( \log{\zeta_{n,k}(1)} \right) \to 0$, $n \to \infty$.
To prove that
$$
\frac{2}{n^2} \sum_{i,j=1}^{n} \cov \left( \log{\phi_{m,l}(i)}, \log{\zeta_{n,k}(j)} \right) \to 0, \,\, n,m \to \infty,
$$
we write, for $i,j=1,\ldots,n$, $u,w > 0$, $x, y \in \mathbb{R}^d, \, x \neq y$, $\norm{x-y} > r_{n-1}(w)$ (thus $n > \frac{w}{\norm{x-y}^d} + 1$) and $m \in \mathbb{N}$,
\vspace{-0.2cm}
\begin{align*}
	\begin{gathered}
		\p\left( \phi_{m,l}(i) \leq u, \zeta_{n,k}(j) \leq w | X_i = x, X_j = y\right) \\ =
		\p\left( \norm{x - Y_{(l)}(x, \mathbb{Y}_m)} \leq r_m(u), \norm{y-X_{(k)}(y, \{X_s\}_{s \neq i,j} \cup \{x\})} \leq r_{n-1}(w) \right) \\ = 1 - \p\left( \norm{x - Y_{(l)}(x, \mathbb{Y}_m)} > r_m(u) \right) - \p\left( \norm{y-X_{(k)}(y, \{X_s\}_{s \neq i,j} \cup \{x\})} > r_{n-1}(w) \right) \\+ \p\left( \norm{x - Y_{(l)}(x, \mathbb{Y}_m)} > r_m(u), \norm{y-X_{(k)}(y, \{X_s\}_{s \neq i,j} \cup \{x\})} > r_{n-1}(w) \right)
 \end{gathered}
\end{align*}
\begin{align*}
	\begin{gathered}
=1 - \sum_{s_1 = 0}^{l-1} \binom{m}{s_1} (W_{m,x}(u))^{s_1} (1 - W_{m,x}(u))^{m-s_1} - \sum_{s_2 = 0}^{k-1} \binom{n-2}{s_2} (V_{n-1,y}(w))^{s_2} (1 - V_{n,y}(w))^{n-2-s_2} \\ +
		\sum_{s_1=0}^{l-1} \sum_{s_2=0}^{k-1} \binom{m}{s_2} (W_{m,x}(u))^{s_1} \binom{n-2}{s_2} (V_{n-1,y}(w))^{s_2} (1-W_{m,x}(u))^{m-s_1} (1-V_{n-1,y}(w))^{n-2-s_2}\\
 =		\left( 1 - \sum_{s_1 = 0}^{l-1} \binom{m}{s_1} (W_{m,x}(u))^{s_1} (1 - W_{m,x}(u))^{m-s_1} \right) \\ \cdot \left( 1 - \sum_{s_2 = 0}^{k-1} \binom{n-2}{s_2} (V_{n-1,y}(w))^{s_2} (1 - V_{n,y}(w))^{n-2-s_2} \right).
	\end{gathered}
\end{align*}
Further we combine the estimates obtained at \text{Steps 4} and \text{5} of Theorem \ref{th_main2} proof. Note that now we  consider
$(x,y) \in A_1 \cap \widetilde{A}_1$ and
employ $G_{\max\{N_1, N_2\}}(\cdot)$.

Thus  we have established that $\var\big(\widehat{D}_{n,m}(k,l)\big) \to 0$ as $n,m\to \infty$, hence \eqref{main2} holds.
The proof is complete. $\square$

\appendix

\section{Proofs of auxiliary results}
Proofs of Lemmas \ref{lemma1}, \ref{lemma_G} and \ref{l1} are similar to the proofs of Lemma 2.5 and 3.1, 3.2 in \cite{Bul_Dim}. We provide them for the sake of completeness.

{\it Proof of Lemma \ref{lemma1}}.

1) Note that $\log\|x-y\|> e_{[N-1]}\geq 1$ if $\|x-y\|>e_{[N]}$ and $N\in \mathbb{N}$.
Hence, for such $x,y$, one has $(\log\|x-y\|)^{\nu}\leq (\log\|x-y\|)^{\nu_0}$ if $\nu    \in (0, \nu_0]$. If $N\geq N_0$ then
$G_N(u)\leq G_{N_0}(u)$ for $u\geq e_{[N-1]}\geq e_{[N_0-1]}$. Thus $K_{p,q}(\nu, N) \leq K_{p,q}(\nu_0, N_0) < \infty$  for $\nu \in (0, \nu_0]$ and any integer $N \geq N_0$.

2) Assume that $Q_{p,q}(\varepsilon_1, R_1) < \infty$. Consider $Q_{p,q}(\varepsilon_1, R)$ where $R>0$.
If $0<R\leq R_1$ then, for each $x\in \mathbb{R}^d$, according to the definition of $M_q$ one has
$M_q(x,R)\leq M_q(x,R_1)$. Consequently, $Q_{p,q}(\varepsilon_1, R)\leq Q_{p,q}(\varepsilon_1, R_1)<\infty$.
Let now $R>R_1$. One has
$$
M_q(x,R)
\leq \max\left\{M_q(x,R_1), \sup_{R_1<r\leq R}\frac{\int_{B(x,R_1)}q(x)dx + \int_{B(x,r)\setminus B(x,R_1)}q(x)dx}{\mu(B(x,r))}\right\}
$$
$$
\leq \max\left\{M_q(x,R_1), M_q(x,R_1)+\frac{1}{\mu(B(x,R_1))}\right\} = M_q(x,R_1)+\frac{1}{\mu(B(x,R_1))}.
$$
Therefore
\begin{gather*}
Q_{p,q}(\varepsilon_1, R) = \int_{\mathbb{R}^d} (M_q(x,R))^{\varepsilon_1} p(x) \, dx \leq \int_{\mathbb{R}^d} \left(M_q(x,R_1)+\frac{1}{R_1^d V_d} \right)^{\varepsilon_1} p(x) \, dx \\ \leq  \max\{1,2^{\varepsilon_1-1}\}\left(Q_{p,q}(\varepsilon_1, R_1) + (R_1^d V_d)^{-\varepsilon_1}\right)< \infty.
\end{gather*}

Suppose now that $Q_{p,q}(\varepsilon_1,R)<\infty$ for some $\varepsilon_1>0$ and $R>0$. Then, for
any $\varepsilon \in (0,\varepsilon_1]$,
the Lyapunov inequality yields $Q_{p,q}(\varepsilon,R)\leq (Q_{p,q}(\varepsilon_1,R))^{\frac{\varepsilon}{\varepsilon_1}}<\infty$.

3)
Let $T_{p,q}(\varepsilon_2, R_2) < \infty$. Take $0<R\leq R_2$. Then, for each $x\in \mathbb{R}^d$, according to the definition of $m_q$ we get
$0\leq m_q(x,R_2) \leq m_q(x,R)$. Hence $T_{p,q}(\varepsilon_2,R)\leq T_{p,q}(\varepsilon_2,R_2)<\infty$.
Consider $R>R_2$. For each $x\in \mathbb{R}^d$ and every $a>0$, the function
$I_q(x,r)$
is continuous in $r$ on $(0,a]$.
Consider an arbitrary (fixed) $x\in S(q)\cap \Lambda(q)$. Then
there exists $\lim_{r\to 0+}I_q(x,r)=q(x)$. For such $x$, set $I_q(x,0):=q(x)$. Thus $I_q(x,\cdot)$ is continuous on any segment $[0,a]$. Hence, one can find $\widetilde{R}_2$ in $[0,R_2]$ such that
$m_q(x,R_2)= I_q(x,\widetilde{R}_2)$ and there exists $R_0$ in $[0,R]$ such that
$m_q(x,R)=I_q(x,R_0)$. If $R_0 \leq R_2$ then $m_q(x,R)=m_q(x,R_2)$ (since $m_q(x,R)\leq m_q(x,R_2)$ for $R>R_2$ and $m_q(x,R)=I_q(x,R_0) \geq m_q(x,R_2)$ as $R_0 \in [0,R_2]$).
Assume that
$R_0 \in (R_2,R]$. Obviously $R_0>0$ as $R_2>0$. One has
\begin{gather*}
m_q(x,R)= I_q(x,R_0)=\frac{\int_{B(x,R_2)}q(y)dy+\int_{B(x,R_0)\setminus B(x,R_2)}q(y)dy}{\mu((x,R_0))} \\
\geq \frac{\int_{B(x,R_2)}q(y)dy}{\mu(B(x,R_0))}= \frac{\mu(B(x,R_2))}{\mu(B(x,R_0))}I_q(x,R_2)\geq \frac{\mu(B(x,R_2))}{\mu(B(x,R_0))}m_q(x,R_2) \\
= \left(\frac{R_2}{R_0}\right)^d m_q(x, R_2) \geq \left(\frac{R_2}{R}\right)^d m_q(x, R_2).
\end{gather*}
Thus in all cases ($R_0 \in [0, R_2]$ and $R_0\in (R_2,R]$) one has $m_q(x,R) \geq \left(\frac{R_2}{R}\right)^d m_q(x,R_2)$ as $R_2 < R$. Taking into account
the relation
$\mu(S(q)\setminus (S(q)\cap \Lambda(q)))=0$ we come to the inequality
$$
T_{p,q}(\varepsilon_2,R)\leq \left(\frac{R}{R_2}\right)^{\varepsilon_2 d}T_{p,q}(\varepsilon_2,R_2)<\infty.
$$
Assume now that $T_{p,q}(\varepsilon_2,R)<\infty$ for some $\varepsilon_2>0$ and $R>0$. Then, for
any $\varepsilon \in (0,\varepsilon_2]$,
the Lyapunov inequality yields $T_{p,q}(\varepsilon,R)\leq (T_{p,q}(\varepsilon_2,R))^{\frac{\varepsilon}{\varepsilon_2}}<\infty$. The proof is complete. $\square$

\vskip0.2cm

{\it Proof of Lemma \ref{lemma_G}}. We start with relation 1).
Note that if a function $g$ is measurable and bounded on a finite interval $(a,b]$ and $\nu$ is a finite measure on the Borel subsets of $(a,b]$ then
$\int_{(a,b]}g(x)\nu(dx)$ is finite. Thus, for each $a\in \left(0,\frac{1}{e_{[N]}}\right]$, using the integration by parts formula (see, e.g., \cite{Shiryaev}, p. 245) we get
\begin{align}\label{lGcor1}
\begin{gathered}
\int_{\left(a, \frac{1}{e_{[N]}}\right]} F(u) \left(-g_N(u)\right) \, du =
\int_{\left(a, \frac{1}{e_{[N]}}\right]} F(u)d \left(-G_N(-\log u)\right) \\
= G_{N}(-\log a) F(a) + \int_{\left(a, \frac{1}{e_{[N]}}\right]} G_N(-\log u) \, d F(u).
\end{gathered}
\end{align}
Assume now that
$\int_{\left(0, \frac{1}{e_{[N]}}\right]} G_N(-\log u) \,dF(u)<\infty$.
Then by the monotone convergence theorem
\begin{equation}\label{c1}
\lim_{a\to 0+}\int_{(0, a]} G_N(-\log u) \, d F(u)=0.
\end{equation}
Clearly, the following nonnegative integral admits an estimate
\begin{gather*}
\int_{(0, a]} G_N(-\log u) \,dF(u) \geq G_N(-\log a) \int_{(0, a]} dF(u)  \\
= G_N(-\log a) (F(a) - F(0)) =
G_N(-\log a) F(a) \geq 0.
\end{gather*}
Therefore \eqref{c1} implies that
\begin{equation}\label{c2}
G_N(-\log a) F(a) \to 0,\;\;a\to 0+.
\end{equation}
Letting $a\to 0+$ in \eqref{lGcor1} we come, by the monotone convergence theorem, to  relation 1) of
our Lemma.
Suppose now that
\begin{equation}\label{c3}
\int_{\left(0, \frac{1}{e_{[N]}}\right]} F(u) \left( -g_N(u) \right) \, du <\infty.
\end{equation}
In view of \eqref{c3} and the equality $\int_{\left(0, \frac{1}{e_{[N]}}\right]} F(u) \left( -g_N(u) \right) \, du = \int_{\left(0, \frac{1}{e_{[N]}}\right]} F(u)d \left(-G_N(-\log u)\right)$ by monotone convergence theorem we have
$
\lim_{b\to 0+}\int_{(0, b]}F(u) \, d(-G_N(-\log u)) = 0.
$
For any $c\in (0,b)$, we obtain the inequalities
$$
\int_{(0, b]}F(u)d(-G_N(-\log u))
\geq \int_{(c, b]}F(u)d(-G_N(-\log u))
$$
$$
=-F(b) G_N(-\log b) + F(c) G_N(-\log c)
+ \int_{(c,b]} G_N(-\log u) \,dF(u)
$$
$$
\geq F(c) G_N(-\log c) - F(b) G_N(-\log b)
+ (F(b)-F(c)) G_N(-\log b)
$$
$$
= F(c) G_N(-\log c) \left(1 - \frac{G_N(-\log b)}{G_N(-\log c)} \right).
$$
Let $c=b^2$ ($b\leq\frac{1}{e_{[N]}} < 1$). Then, for all positive $b$ small enough,
\begin{align*}
    \begin{gathered}
        1 - \frac{G_N(-\log b)}{G_N(-\log c)} = 1 - \frac{G_N(-\log b)}{G_N(-2 \log b)} = 1 - \left(\frac{1}{2}\right) \frac{\log_{[N]}(-\log b)}{\log_{[N]}(-2\log b)} \geq \frac{1}{2}.
    \end{gathered}
\end{align*}
Thus $\int_{(0, b]}F(u)d(-G_N(-\log u))
\geq \frac{1}{2} F(b^2) G_N(-\log (b^2)) \geq 0$.
It follows that
$F(b^2) G_N(-\log b^2) \to 0$ as $b\to 0$. Hence we come to \eqref{c2} taking $a=b^2$.
Then \eqref{lGcor1} yields relation 1).

If one of (nonnegative) integrals appearing in 1) is infinite and other one is finite
we come to the contradiction. Hence 1) is established.
In a similar way one can prove that relation 2) is valid. Therefore, we omit further details.
$\square$

{\it Proof of Lemma \ref{l1}}.
Take $x\in S(q)\cap \Lambda(q)$ and $R>0$. Assume that $m_q(x,R)=0$.
Since the function $I_q(x,r)$ defined in \eqref{I} is continuous in $(x,r)\in \mathbb{R}^d\times (0,\infty)$, there exists $\widetilde{R}\in [0,R]$
($\widetilde{R}=\widetilde{R}(x,R)$) such that $m_q(x,R)= I_q(x,\widetilde{R})$ (recall that  $I_q(x, 0) := \lim_{r \rightarrow 0+} I_q(x,r) = q(x)$ for all $x \in \Lambda(q)$ by continuity).
If $\widetilde{R}=0$ then $m_q(x,r)=q(x)>0$ as $x\in S(q)\cap \Lambda(q)$.
Hence we have to consider $\widetilde{R} \in (0,R]$.
If $I_q(x,\widetilde{R})=0$ then $\int_{B(x,r)}q(y)dy =0$ for any $0<r\leq \widetilde{R}$.
Thus \eqref{2} ensures that $q(x)=0$. However, $x\in S(q)\cap \Lambda(q)$.
So $m_q(x,R) > 0$ for $x \in S(q) \cap \Lambda(q)$. Thus, $S(q)\cap \Lambda(q)\subset D_q(R):=\{x\in S(q): m_q(x,R)>0\}$. It remains to note
that $S(q)\setminus \Lambda(q)\subset \mathbb{R}^d \setminus \Lambda(q)$ and $\mu(\mathbb{R}^d \setminus \Lambda(q))=0$. Therefore $\mu(S(q) \setminus D_q(R)) = 0$. $\square$

\vskip0.2cm

{\it Proof of Lemma \ref{l4}}.
We verify that, for given $N\in \mathbb{N}$ and $\tau >0$,  there exist $a:=a(\tau)\geq 0$ and $b:=b(N,\tau)\geq 0$ such that, for any $c \geq 0$,
\begin{equation}\label{eq_123}
G_N(\tau c) \leq a G_N(c) + b.
\end{equation}
For $c=0$ the statement is obviously true. Let $c>0$.
One can easily see that $\frac{\log_{[N]}(\tau c)}{\log_{[N]}(c)} \to 1$ as $c \to \infty$.
Hence one can find $c_0(N, \tau)$ such that, for all $c \geq c_0(N,\tau)$, the inequality $\frac{\log_{[N]}(\tau c)}{\log_{[N]}(c)} \leq 2$ is valid. Consequently, for  $c \geq c_0(N,\tau)$,
$$
    \frac{G_N(\tau c)}{G_N(c)} = \frac{\tau c \log_{[N]}(\tau c)}{c \log_{[N]}(c)} \leq 2 \tau := a(\tau).
$$
For all $0 \leq c \leq c_0(N, \tau)$ we write $G_N(\tau c) \leq G_N (\tau c_0(N,\tau)) := b(N,\tau)$.
Therefore, for any $c \geq 0$, we come to \eqref{eq_123}.
Thus, for any $\nu > 0$ and $x,y\in \mathbb{R}^d$, $x\neq y$, one has $$
G_{N}(|\log(\|x-y\|^d)|^{\nu}) = G_{N}(d^{\nu}|\log(\|x-y\|)|^{\nu}) \leq a(d^{\nu}) G_{N}(|\log(\|x-y\|)|^{\nu}) + b(N, d^{\nu}).\;\;\square
$$

\vskip0.2cm

{\it Proof of Lemma \ref{l6}}.
For $t \in [0,e_{[N-1]}]$, a function $G_N(t) \equiv 0$ is convex.
We show that $G_N$ is convex on $(e_{[N-1]},\infty)$.
Consider $t > e_{[N-1]}$. Write $\prod\limits_{\varnothing} := 1$ and $\sum\limits_{\varnothing} := 0$. Then, for $N\in \mathbb{N}$,
\vspace{-0.4cm}
\begin{gather*}
    (G_{N}(t))' = \log_{[N]}(t) + \prod_{j=1}^{N-1} \frac{1}{\log_{[j]}(t)}.
\end{gather*}
Obviously,  $\left(\frac{1}{\log_{[k]}(t)}\right)' = -\frac{1}{t\log_{[k]}^2(t)} \prod_{s=1}^{k-1} \frac{1}{\log_{[s]}(t)}$, $k\in \mathbb{N}$. Thus, for $t>  e_{[N-1]}$, we get
\begin{gather*}
    \left(G_N(t)\right)'' = \frac{1}{t} \prod_{j=1}^{N-1} \frac{1}{\log_{[j]}(t)} + \sum_{k=1}^{N-1} \left( -\frac{1}{t} \frac{1}{\log_{[k]}^2 (t)} \prod_{s=1}^{k-1} \frac{1}{\log_{[s]}(t)} \prod_{j\in \{1,\ldots,N-1\}, j \neq k} \frac{1}{\log_{[j]}(t)} \right)\\ =
    \frac{1}{t} \left( \prod_{j=1}^{N-1} \frac{1}{\log_{[j]}(t)} \right)  \left( 1-\sum_{k=1}^{N-1} \prod_{s=1}^k \frac{1}{\log_{[s]}(t)} \right).
\end{gather*}
For $N=1$ and $t>0$, we have $\left(G_1(t)\right)'' = \frac{1}{t} > 0$. Take now $N > 1$.
Clearly, for $t > e_{[N-1]}$, one has $\frac{1}{t} \prod\limits_{j=1}^{N-1} \frac{1}{\log_{[j]}(t)} > 0$  because $\log_{[j]}(t) > \log_{[j]}(e_{[N-1]}) = e_{[N-1-j]} \geq 1 > 0$ when $1 \leq j \leq N-1$.
Observe also that
\begin{equation}\label{last}
    \sum_{k=1}^{N-1} \prod_{s=1}^{k} \frac{1}{\log_{[s]}(t)} < \sum_{k=1}^{N-1} \prod_{s=1}^{k} \frac{1}{e_{[N-1-s]}} \leq \sum_{k=1}^{N-1} \frac{1}{e_{[N-2]}} = \frac{N-1}{e_{[N-2]}} \leq 1.
\end{equation}
The last inequality is established by induction in $N$.
Thus, in view of \eqref{last}, we have proved that, for all  $t > e_{[N-1]}$ and $N \in \mathbb{N}$, the inequality $(G_N(t))'' > 0$ holds. Hence, the function $G_N(t)$ is (strictly) convex on $\left(e_{[N-1]}, \infty\right)$.

Let $h:[a,\infty)\to \mathbb{R}$ be a continuous nondecreasing function. If the restrictions of $h$ to $[a,b]$ and $(b,\infty)$ (where
$a<b$) are convex functions then, in general, it is not true that $h$ is convex on $[a,\infty)$.
However, we can show that $G_N$ is convex on $[0,\infty)$. Note that a function $G_N$ is convex on
$[e_{[N-1]},\infty)$ since it is convex on $(e_{[N-1]},\infty)$ and continuous on $[e_{[N-1]},\infty)$. Take now any $z \in [0, e_{[N-1]}]$, $y \in (e_{[N-1]}, \infty)$ and $s \in [0,1]$. Then $G_N(s z + (1-s) y) \leq G_N(s e_{[N-1]} + (1-s) y) \leq s G_N(e_{[N-1]}) + (1-s) G_N(y) = (1-s) G_N(y) = s G_N(z) + (1-s) G_N(y)$ as $G_{N}(z) = 0$.
Thus, for each $N\in\mathbb{N}$, a function $G_N(\cdot)$ is convex on $\mathbb{R}_{+}$. $\square$

\vskip0.2cm
{\it Proof of Corollary \ref{cor4}}.
The proof (i.e. checking the conditions of both Theorem \ref{th1} and \ref{th_main2}) is quite similar to the proof of Corollary 2.11 in \cite{Bul_Dim}.

\vspace{0.2cm}
{\bf Acknowledgements}
The authors are grateful to Professor A.Tsybakov for useful discussions.
This work is supported by the Lomonosov Moscow State University under grant ``Modern Problems of the Fundamental Mathematics and Mechanics''.

%
%
%
%

\end{document}